\documentclass[aap,preprint]{imsart}
\pdfoutput=1
\setattribute{journal}{name}{}
\arxiv{math.PR}
\usepackage{amssymb,enumerate}
\usepackage{amsfonts}
\usepackage{natbib}

\usepackage{amsmath,amsthm}

\newcommand{\field}[1]{\mathbb{#1}}
\DeclareMathOperator{\E}{\field{E}}              
\def\R{\field{R}}                                
\def\supp{\text{supp}}

\newtheorem{thm}{Theorem}[section]
\newtheorem{cor}[thm]{Corollary}
\newtheorem{lem}[thm]{Lemma}
\newtheorem{prop}[thm]{Proposition}
\newtheorem{ex}[thm]{Example}
\newtheorem{defn}[thm]{Definition}
\newtheorem{rem}[thm]{Remark}

\def\cnvg{\stackrel{v}{\to}}

\def\bZ{\boldsymbol Z}

\newcommand{\nwc}{\newcommand}
\nwc{\COM}[1]{}

\def\bC{\mathbb{C}}
\def\C{\mathcal{C}}

\def\bO{\mathbb{O}}
\def\O{\mathcal{O}}

\def\bone{\boldsymbol 1}
\def\ba{\boldsymbol a}

\def\bzero{\boldsymbol 0}

\def\bx{\boldsymbol x}
\def\bz{\boldsymbol z}

\def\by{\boldsymbol y}
\def\binfty{\boldsymbol \infty}

\def\bE{\mathbb{E}}

\def\bM{\mathbb{M}}
\def\R{\mathbb{R}}
\def\C{\mathcal{C}}
\def\D{\mathbb{D}}

\def\bF{\mathbb F}

\def\bC{\mathbb{C}}
\def\bD{\mathbb{D}}

\begin{document}

\begin{frontmatter}

\title{Living on the multi-dimensional edge: Seeking hidden risks using regular variation}

\runtitle{Hidden risks}

\begin{aug}
  \author{\fnms{Bikramjit}  \snm{Das}\thanksref{m1}\ead[label=e1]{bikram@math.ethz.ch}},
  \author{\fnms{Abhimanyu}  \snm{Mitra}\thanksref{m2}\ead[label=e2]{am492@cornell.edu}}
\and
  \author{\fnms{Sidney}  \snm{Resnick}\thanksref{m2}\ead[label=e3]{sir1@cornell.edu}}

  \thankstext{m1}{B. Das was partially supported by the IRTG, Zurich-Berlin on {\it{Stochastic Models for Complex Processes}}.}
  \thankstext{m2}{S. Resnick and A. Mitra were partially supported by ARO Contract W911NF-10-1-0289 at Cornell University.}

  \runauthor{Das, Mitra \and Resnick}

  \affiliation{ETH Zurich, Cornell University \and Cornell University}

  \address{RiskLab, Department of Mathematics,\\ETH Zurich, \\R\"amistrasse 101,\\ 8092 Zurich, Switzerland\\
           \printead{e1}}

  \address{School of Operations Research and\\ Information Engineering,\\ Cornell University,\\
Ithaca, NY 14853, USA\\
          \printead{e2,e3}}

\end{aug}

\begin{abstract}  {Multivariate regular variation plays a  role  assessing tail risk
    in diverse applications such as finance, telecommunications,
    insurance and environmental science. The classical theory, being
    based on an asymptotic model,
sometimes leads to inaccurate and useless estimates of probabilities of joint tail
    regions. This problem can be  partly ameliorated by using
    {\it  hidden
    regular variation\/} \citep{resnick:2002a, mitra:resnick:2010}.
We offer a more flexible definition of hidden regular variation
  that provides improved risk estimates for  a larger class
    of risk tail regions.}
\end{abstract}



\end{frontmatter}

\section{Introduction}\label{sec:intro}
Daily we observe  environmental, technological and
  financial phenomena possessing inherent risks.
There are
financial risks {from large investment} losses;
environmental risks from health hazards resulting from high 
concentrations of atmospheric pollutants; hydrological risks from
river floods.   {Risk analysis requires estimation of}
tail probabilities that provide measures of such risks.  The
mathematical framework of multivariate regular variation provides 
tools to compute tail probabilities associated with such
risks; see \cite{resnickbook:2007, joe:li:2010, cai:einmahl:dehaan:2011}. {These tools have
  limitations which we begin to address in this paper.}

Consider a non-negative random vector
${\bZ } = (Z^1, Z^2, \cdots, Z^d)$ called a {\it risk vector\/}. The
distribution of ${\bZ }$ has multivariate regular variation if there exist a function
$b(t) \uparrow \infty$ and a non-negative non-degenerate Radon measure
$\mu(\cdot)$ on $\E = [0, \infty]^d \setminus \{ (0, 0, \cdots, 0) \}$
such that as $t \to \infty$,
\begin{equation}\label{eqn:existing_reg_var}
tP\left[ \frac{{\bZ }}{b(t)} \in \cdot \right] \stackrel{v}{\rightarrow} \mu(\cdot),
\end{equation} 
 where $\stackrel{v}{\rightarrow}$ denotes vague convergence in
 $\bM_+(\E)$,
the set of all Radon measures on $\E$ \citep[page
172]{resnickbook:2007}.
Note  that \eqref{eqn:existing_reg_var} effectively assumes tail equivalence
of the marginal components \citep[Section
6.5.6]{resnickbook:2007}, so while  \eqref{eqn:existing_reg_var} is
valuable as a theoretical foundation it must be modified for applications.

The asymptotic relation \eqref{eqn:existing_reg_var} allows
 the limit measure $\mu(\cdot)$  to be
 used for approximating tail probabilities. For example, 
{approximation of}
the probability 
of the event $\{Z^i>x^i \text{ for some } i \}$ for large thresholds
$x^i, i = 1, 2, \cdots, d$ requires {the} ability to  compute $\mu\left(
  \left \{ (z^1, \cdots, z^d) \in \E : z^i > w^i \text{ for some } i
  \right\} \right)$ for $w^i > 0$, $i = 1, 2, \cdots, d$.  
Such approximations of tail probabilities are  sensitive to
degeneracies in the limit measure $\mu(\cdot)$. 
 For example, {when {\it asymptotic independence\/} is present as in Gaussian
 copula models,}  the limit measure $\mu(\cdot)$ in
\eqref{eqn:existing_reg_var} concentrates on the coordinate axes 
$\mathbb{L}_i := \{\bx\in \R^d: x^j=0 \; \forall j \neq i \},\,
i=1,\dots ,d,$ and 
$\mu\left(\{(z^1,z^2, \cdots, z^d) \in \E : z^{i} > w^1, z^{j} > w^2 \right) = 0$
for any $1 {\le} i < j {\le} d$ and $w^1,w^2>0$. Consequently, we {would}
approximate the joint tail probability  
 $$ P(Z^{i} > x^1, Z^{j} > x^2) \approx 0 $$ 
for large thresholds $x^1,x^2$ and {conclude risk contagion is absent. }
This {conclusion may be} naive and hence the concept of
{\it{hidden regular variation}} (HRV) was introduced 
\citep{resnick:2002a}
 which offered a refinement of this approximation; see
 \cite{maulik:resnick:2005, heffernan:resnick:2005,
   mitra:resnick:2010} and the seminal concept of {\it{coefficient of
     tail dependence}}  in \cite{ledford:tawn:1996,
   ledford:tawn:1998}.

{The definition of hidden regular variation} offers some strengths
but also has weaknesses.
{The existing definition} provides insight only in the presence of a
restricted class 
 of degeneracies in the limit measure $\mu(\cdot)$ in \eqref{eqn:existing_reg_var};
namely when $\mu(\cdot)$  concentrates either on the coordinate axes,
or the coordinate planes or similar coordinate hyperplanes in higher
dimensions.  {However, other degeneracies in $\mu(\cdot)$ are possible; for example,
$\mu(\cdot)$ may be concentrated on the diagonal $\{ (z^1, z^2,
\cdots, z^d) \in \E : z^1 = z^2 = \cdots = z^d \}$, {a condition
  called asymptotic full dependence}. To deal with such degeneracies
and situations where $\mu $ may place zero mass on large portions of
the state space,
we define  in Section \ref{sec:rvandhrv}
hidden regular variation
on cones.  For us, a cone $\mathbb{C}$ in $\R^d$ is a set
$\mathbb{C} \subset \R^d$ satisfying ${\bf{x}} \in \mathbb{C}$
implies $t{\bf{x}} \in \mathbb{C}$ for $t > 0$. }

In practice, different risk assessment 
problems {require} calculating tail probabilities for different
kinds of events. {Hidden regular variation, as previously defined,
may be a natural choice for some calculations but not for others.}
For example,  suppose  $(Z^1,Z^2)\in \R_+^2$
is a risk vector and we must calculate
risk probabilities of the form $P\left(| Z^1 - Z^ 2| > w \right)$ for large thresholds $w > 0$. 
{I}f we use multivariate
regular variation when 
the limit measure $\mu(\cdot)$ in \eqref{eqn:existing_reg_var} 
is concentrated on the diagonal $\left \{ (z^1, z^2) \in \E: z^1 = z^2
\right \}$, the tail probability $P\left(| Z^1 - Z^ 2| > w \right)$
{must be} approximated as zero {for large thresholds $w > 0$.}
{The existing notion of hidden regular variation designed to help when
$\mu $ concentrates on the axes}
cannot offer a refinement in this case.
Example \ref{ex:asym_full_dep} of Section \ref{sec:examples}
{illustrates} how a more general  theory  overcomes this difficulty.  
  This along with other examples in Section \ref{sec:examples}
  emphasize the need for the theory of hidden regular variation on
  general cones.

The conditional extreme value (CEV) model
\citep{heffernan:tawn:2004, heffernan:resnick:2007, das:resnick:2011}
provides {one} alternative approach to multivariate extreme value
modeling.  In  standard form, the CEV model can also be formulated as
 regular variation on a particular cone in  $\E$ and this is discussed in Section
\ref{sec:CEVandNon-stanMRV_connection}. Also in   Section
\ref{sec:CEVandNon-stanMRV_connection}, we consider
non-standard regular variation 
from the point of view of regular variation on a sequence of cones.
Non-standard regular variation is essential in practice
since in applications we cannot assume tail-equivalence of all
marginals in a multivariate model as is done in
\eqref{eqn:existing_reg_var}.

Section \ref{sec:detection_estimation} discusses how to fit our model of HRV on cones to data as well as estimation 
techniques of tail probabilities using our model. {We have 
adapted ideas previously used in multivariate heavy tail analysis and
this discussion is not comprehensive, merely a feasibility display. In
particular we have not performed data analyses.}
We close our
discussion with concluding remarks in Section
\ref{sec:conclusion} and give some deferred results and proofs in 
Section \ref{sec:appendix}.

Our definition of HRV relies on a notion of convergence of
measures called
$\mathbb{M}^*$-convergence
that is similar to  $\mathbb{M}_0$ convergence of
\cite{hult:lindskog:2006a}. This
$\mathbb{M}^*$-convergence is developed in Section
\ref{sec:conv_concepts} {where we  also discuss  reasons to
abandon} the standard practice of defining regular variation
through vague convergence on a  compactification of $\R^d$.

\subsection{Notation}\label{notations} We briefly discuss
some frequently used notation
and concepts.

\subsubsection{Vectors, norms and topology} Bold letters are used
to denote vectors, with capital letters reserved for random vectors
and small letters for non-random vectors, e.g.,
 $\bx = (x^1, x^2, \cdots, x^d) \in \R^d$. We also denote
$${\bzero }= (0, 0, \cdots, 0), \quad \bone  = (1, 1, \cdots, 1),
\quad {\boldsymbol{\infty}} = (\infty, \infty, \cdots, \infty).$$ 
Operations on and between vectors are understood componentwise. For
example, for vectors $\bx, \bf{z}$, 
\[\bx \, {\le} \, \bz  \quad \text{means} \quad x^i \le z^i, \; i=1,\ldots,d.\] 
For a set $A \subset {[0,\infty)}^d$ and $\bx \in A$, 
use $[ \bzero, \bx ]^c$ to mean 
$ [\bzero, \bx ]^c = A \backslash [\bzero, \bx] = \{ \by  \in A:
\vee_{i = 1}^d {y^i}/{x^i} > 1 \}$.   When we use the notation $
[\bzero, \bx ]^c$, the set $A$ should be clear from
the context. 

For the $i$-th largest component of $\bx$, we write $x^{(i)}$, that is,
$x^{(1)} \ge x^{(2)} \ge \cdots \ge x^{(d)}.$  
Thus a superscripted number $i$ denotes the $i$-th component of a
vector, whereas a superscripted $(i)$  denotes the $i$-th largest
component in the vector.

Operations with $\infty$ are understood using the conventions:
\begin{align*}
 & \infty + \infty = \infty, && \infty - \infty = 0, && \text{for } x\in \R, \infty+x=\infty -x = \infty,\\
 & 0. \infty = 0, && \text{for } x>0, \; x.\infty = \infty, &&\text{for } x<0,\;  x.\infty = -\infty.
\end{align*}

{Fix a norm on $\mathbb{R}^d$ and denote the norm of $\bx$ as $|| \bx||$. }
Let  $d(\bx,\by)=\|\bx-\by\|$ be  the metric induced by the norm and, 
as usual, for $A \subset \mathbb{R}^d$, set $d(\bx,A)=\inf_{\by \in A}
d(\bx,\by)$. When attention is focused on the set $\bC$,  and  $A \subset
\bC$, 
the $\delta $-dilation or swelling of $A$ in $\bC$ is 
$ A^\delta := \{x\in \bC: d(x,A) < \delta \}$. The topology on
$\mathbb{R}^d$ is the usual norm topology referred to as the Euclidean topology and the topology on a
{subset of}
$\mathbb{R}^d$  is the relative topology induced by the Euclidean topology.

{Two sets $A$ and $B$ in $\mathbb{R}^d$ are {\it bounded away\/} from each
other if $\bar A \cap \bar B  = \emptyset$, where $\bar
A$ and $\bar B$  
are the closures of $A$ and $B$.}

\subsubsection{Cones}\label{subsub:cones}
We  denote by
$\E=[0,\infty]^d\setminus \{\bzero\}$ and
$\mathbb{D}=[0,\infty)^d\setminus \{\bzero\}$, the one
point puncturing of the compactified and uncompactified
versions of $\R_+^d$. The symbols $\E$ and $\D$ {may} appear with
superscripts 
 denoting subsets of the compactified and uncompactified $\R^d_{+}$
 respectively. For example, $\E^{(l)} =\{\bx \in [0,\infty]^d:
 x^{(l)} > 0\}$ 
 and $\mathbb{D}^{(l)} =\{\bx \in [0,\infty)^d: x^{(l)} > 0\}$,
 where $x^{(l)}$ is the $l$-th largest component of $\bx$. {Note $\E^{(1)}=\E$.}

A set $\bC\subset\R^d$ is a cone if 
$\bx\in \bC$ implies $t\bx \in \bC$ for all $t>0$. Cones in the
Euclidean space are usually denoted by  mathematical bold symbols $\bC,
\bD, \bE, \bF$, etc. Since  one typically deals with non-negative risk
vectors, we focus on the case where  $\bC \subset  [0, \infty)^d$.  
Call a subset $\bF \subset \bC$ a closed cone in $\bC$ if $\bF$ is a
closed subset of $\bC$ as well as a cone. Example: $\bF =\{\bzero\}$
or when $d=2$,  $\bF=\{(t,0): t\geq 0\}$.
{The complement of the closed cone $
\bF$ in { $\bC$} is  an open cone in $\bC$; that is, the complement {of
$\bF$} is
an open subset in $\bC$ as well as a cone.} 

 Fix a closed cone $\bC \subset [0,\infty)^d $  containing $\bzero$
 and suppose $\bF \subset \bC$ be a closed cone in $\bC$ containing
 $\bzero$. 
Then $\bO:=\bC \setminus \bF$  is an open cone and 
 $\bC$ and $\bF$ are complete separable metric spaces under the metric $d(\cdot,\cdot)$.
 Let $\C$ denote the Borel $\sigma$-algebra of $\bC$. 
 Clearly $\bO$ is again a {separable metric space (not necessarily
   complete)} equipped with the $\sigma$-algebra $\O=\{B\subset \bO: B
 \in \C \}$. 
  
Examples: (i) $\bC=[0,\infty)^d$, $\bF=\{\bzero\}$ and then
$\bO=\D=[0,\infty)^d \setminus \{\bzero\}$. (ii) $d=2$,
$\bC=[0,\infty)^2$, $\bF
=\{(0,x): x\geq 0\} \cup \{(y,0): y\geq 0\}$,
$\bO=(0,\infty)^2$. (iii)
 $d=2$,
$\bC=[0,\infty)^2$, $\bF
=\{(x,0): x\geq 0\},$ and $\bO=
\mathbb{D}_{\sqcap}:=[0,\infty)\times(0,\infty)$.  


\subsubsection{Regularly varying
  functions}\label{subsubsec:regVarFcts}
A function $U:[0,\infty)\mapsto [0,\infty)$ is regularly varying with
index $\beta \in \mathbb{R}$ if for all $x>0$,
$$\lim_{t\to\infty} \frac{U(tx)}{U(t)}=x^\beta.$$
We write $U \in RV_\beta$. See \citet{resnickbook:2008,
  dehaan:ferreira:2006,bingham:goldie:teugels:1987}. 

\subsubsection{{{{Vague convergence of measures.}}}}
\label{subsec:rv_vague}
 We express vague convergence of Radon measures as
$\stackrel{v}{\rightarrow}$ (\cite[page 173]{resnickbook:2007},
\cite{kallenberg:1983})
and
weak convergence of probability measures as $\Rightarrow$ \cite[page
14]{billingsley:1999}. Denote the set of non-negative Radon
measures on a space $\mathbb{S}$ as $\bM_+(\mathbb{S})$ and the set of all non-negative
continuous functions with compact support from $\mathbb{S}$ to {$\R_+$} as
$C_K^+(\mathbb{S})$. 

Vague convergence on $\E$  has traditionally been used when defining
multivariate regular variation. We explain in the next section why
continuing with this practice is problematic and what should be done.

\section{$\mathbb{M}^*$-convergence of measures and regular variation on cones} \label{sec:conv_concepts}
  
{This paper {requires} a definition of multivariate regular variation {on}
cones of the Euclidean space which differs from the traditional
definition through {\it{vague convergence}} of measures.  {Following \cite{hult:lindskog:2006a},
we define
regular variation based on a notion of convergence of measures
we call $\mathbb{M}^*$-convergence. }

\subsection{Problems with compactification of $\R^d$} 
Multivariate regular variation on $[0,
\infty)^d$ is usually defined using vague convergence of Radon measures
on $\E=[0, \infty]^d\setminus \{\bzero\}$ \citep{resnickbook:2007}.
{The reason for compactifying $[0,\infty)^d$ and then removing $\bzero$
is this makes sets bounded away from $\{\bzero\}$  relatively compact
(cf. \cite[Section 6.1.3]{resnickbook:2007}) and since 
Radon measures put finite mass on relatively compact sets,
this theory is suitable for estimating probabilities of
tail regions.}

{{The} theory of hidden regular variation {may} require removal of more than just a point.
Furthermore, compactifying from $[0,\infty)^d $ to $[0,\infty]^d $
introduces problems.  For one thing, it is customary to rely
heavily on the polar coordinate transform
$$\bx \mapsto \Bigl(\|\bx\|, \frac{\bx}{\|\bx\|} \Bigr)$$
which is  only defined on $[0,\infty)^d \setminus
\{\bzero\}$ and if the state space $[0,\infty]^d \setminus \{\bzero\}$
is used  an awkward kluge \citep[page 
176]{resnickbook:2007} is required to show the equivalence of regular variation in
polar and Cartesian coordinates. {A} workaround
 is only possible because the limit measure $\mu $ in
\eqref{eqn:existing_reg_var}  puts zero mass on lines through
  infinity $\{\bx: \vee_{i=1}^d x^i =\infty\}$ but absence of mass on
  lines through $\binfty$ does 
  not necessarily persist for regular variation on other cones \citep{mitra:resnick:2010}.}

 {Also, compactification introduces counterintuitive {geometric} properties.
{For example,}
the topology on
  $[0, \infty]^d$ can be defined through a homeomorphic map  $[0,
  \infty]^d \mapsto [0, 1]^d$, {such as}
$$\bz  = (z^1, z^2, \cdots, z^d) \mapsto \bigl(
z^1/(1+z^1) , \dots, z^d/(1+z^d) \bigr) .$$
{Restrict attention to $d=2$} 
and consider  two parallel lines {in} $[0, \infty]^2$ with the same
positive
and finite 
slope. These lines both converge to the same point $(\infty, \infty)$
and therefore in the compactified space,
these two parallel lines are not
 bounded away from each other.  
Interestingly, this is not the case if the lines are horizontal or vertical.

To see the impact that parallel lines not being  bounded away from each
other can have recall one of {the motivational examples from Section \ref{sec:intro} with} $d=2$, where the limit measure
$\mu(\cdot)$ in \eqref{eqn:existing_reg_var} is concentrated on the
diagonal $\mathbb{DIAG} := \{(z^1, z^2) \in \E: z^1 = z^2 \}$ and we
need to approximate the tail probability $P\left(| Z^1 -
  Z^ 2| > w \right)$ for a large threshold $w > 0$. Of course, if
we use multivariate regular variation as in
\eqref{eqn:existing_reg_var} to approximate $P\left(| Z^1 - Z^ 2| > w
\right)$, we approximate $P\left(| Z^1 - Z^ 2| > w \right)$ as
zero. {If $P[Z^1=Z^2]<1$,} this approximation is crude.
Following the usual definition of HRV, 
we remove the diagonal $\mathbb{DIAG}$ and define
 regular variation on the sub-cone
${\left(\mathbb{DIAG}\right)}^c := \{(z^1, z^2) \in \E: z^1 \ne z^2
\}$. {Since we seek}  to {approximate} $P\left(| Z^1 - Z^ 2| > w
\right)$, we
{are interested in the set}  $A_{>w} := \{(z^1, z^2) \in \E:
|z^1- z^2| > w \}$.  If we  define HRV
 on the sub-cone ${\left(\mathbb{DIAG}\right)}^c$ {as an asymptotic
 property using}
vague convergence, we {need} the set $A_{>w}$ to be
relatively compact in the sub-cone
${\left(\mathbb{DIAG}\right)}^c$. However, if the sub-cone
${\left(\mathbb{DIAG}\right)}^c$ is endowed with the relative topology
from the topology on $[0, \infty]^2$,
$A_{>w}$ {is} not relatively compact {since}
 the boundaries of  $A_{>w}$ are the two parallel lines $\{(z^1, z^2) \in \E:
z^1- z^2 = w\}$ and $\{(z^1, z^2) \in \E: z^1- z^2 = - w \}$, which
are both parallel to the diagonal $\mathbb{DIAG}$.
In the topology of $[0, \infty]^d$, the boundaries of the set
$A_{>w}$ are not bounded away from the diagonal $\mathbb{DIAG}$ and
hence by Proposition 6.1 of \cite[page 171]{resnickbook:2007}, the set
$A_{>w}$ is not relatively compact in
${\left(\mathbb{DIAG}\right)}^c$.

As already observed,  horizontal or vertical parallel
lines are bounded away from each other in $[0, \infty]^2$. 
If the limit measure $\mu(\cdot)$ in \eqref{eqn:existing_reg_var}  
concentrates on the axes, the traditional definition of HRV  \citep{resnick:2002a}
removes the axes and defines
hidden regular variation on the cone $(0,\infty]^2$. However, risk
regions of interest of the form $(z^1,\infty]\times (z^2,\infty]$ are
still relatively compact and we do not encounter a problem as
above..

{Thus, we conclude
that a flexible theory of hidden regular variation on general cones of
$[0, \infty)^d$ requires considering the possibility that compactification and vague
convergence be abandoned.}
However, without
compactification, how do we guarantee risk sets corresponding to tail
events are relatively compact
and their probabilities approximable by asymptotic methods?  A theory 
based on  $\mathbb{M}^*$-convergence of measures sidesteps many difficulties.

 \subsection{{$\mathbb{M}^*$}-convergence of measures} \label{sec:Mstar}
We follow ideas of  \cite{hult:lindskog:2006a}
{who removed only}
a fixed point {from a closed set}, whereas we remove a closed cone.  

{As in Section \ref{subsub:cones}
we fix  a closed cone $\bC \subset [0,\infty)^d$ containing $\bzero$ 
and $\bF \subset \bC$ is a closed cone in $\bC$ containing $\bzero$.
Set $\bO:=\bC \setminus \bF$, which is an open cone in $\bC$.
 Let $\C$ {be} the Borel $\sigma$-algebra of $\bC$ and  the
 $\sigma$-algebra in 
 $\bO$ is $\O=\{B\subset \bO: B \in \C \}$. }  
Denote by $\C_{\bF}$ the  set of all bounded, continuous real-valued functions $f$ on $\bC$ such that
 $f$  vanishes on 
$\bF^r:= \{\bx \in \bC: d(\bx,\bF)<r\}$
for some $r>0$.
The class of Borel measures on $\O$ that assign finite measure to all $D \in \O$
which are bounded away {from} $\bF$   {is called}
$\bM^*(\bC,\bO)$.
  Equivalently, $\mu \in \bM^*(\bC,\bO)$ if and only if
$\mu$ is finite on $\bC \setminus \bF^r  $ for all $r>0$.

\begin{defn}[$\bM^*(\bC,\bO)$-convergence]\label{def:M*conv}
{\rm
 For $\mu,\mu_n\in  \bM^*(\bC,\bO), n\ge 1$, $\mu_n$ converges to $\mu$ in $\bM^*(\bC,\bO)$ if 
\begin{align}\label{conv:M^*}
 \lim_{n \to \infty} \mu_n(B) = \mu(B),
\end{align}
for all $B \in {\mathcal{O}}$ with $\mu(\partial B) = 0$ and 
{$B$ bounded away from $\bF$.}
 We write $\mu_n \stackrel{*}{\rightarrow} \mu$ in $\bM^*(\bC,\bO)$ as $n\to \infty$.}
\end{defn}

We can  {metrize} the space $\bM^*(\bC,\bO)$. One method: For $\mu, \nu \in
\bM^*(\bC,\bO)$ define 
\begin{align}\label{def:metric}
 d_{\bM^*}(\mu,\nu)&:= \int_0^{\infty} e^{-r} \frac{d_P(\mu^{(r)}, \nu^{(r)})}{1+d_P(\mu^{(r)}, \nu^{(r)})} dr,
\end{align}
where 
{$\mu^{(r)},\nu^{(r)}$ are  the
restrictions of $\mu, \nu$  to  $\bC\setminus \bF^r $ and $d_P$ is the
Prohorov metric \citep{prohorov:1956}. }

Following \cite{hult:lindskog:2006a}, 
$(\bM^*(\bC,\bO),d_{\bM^*})$  is a complete separable metric space
and the 
expected analogue of the Portmanteau theorem
\citep{billingsley:1999} holds:
For $\mu_n \in  \bM^*(\bC,\bO), n\ge 0$, the following are equivalent:
\begin{enumerate}
\item  $\mu_n \stackrel{*}{\rightarrow} \mu_0$ in $\bM^*(\bC,\bO)$ as $n \to \infty$. 
\item For each $f \in \C_{\bF}$, $\int f d \mu_n \to \int f d\mu_0$ as $n \to \infty$.
\item $\limsup_{n\to \infty} \mu_n(A) \le \mu_0(A)$ and $\liminf_{n\to
    \infty} \mu_n(G) \ge \mu_0(G)$ for all closed sets $A \in \O$ 
and  open sets $G \in \O$ such that $ \overline{G} \cap \bF =\emptyset$.
\end{enumerate}

 \section{Regular and hidden regular variation on
   cones} \label{sec:rvandhrv} 
We define regular variation on a nested sequence of cones, where each
cone is a subset of the previous one. Each cone in the sequence
possesses a different regular variation, which remains hidden while
studying regular variation on the bigger cones in the sequence.

\subsection{Regular variation}
We use the concepts of Section \ref{sec:conv_concepts} to define
regular variation.

\begin{defn}\label{def:closedconecomp}
 {\rm
{Suppose $\bF \subset \bC \subset [0,\infty)^d$ and $\bF$ and
       $\bC$ are closed cones containing $\bzero$.
  A random vector ${\bZ } \in \bC$ {has a distribution with a} regularly 
varying {tail} on $\bO = \bC \setminus \bF$, 
 if there {exist} a function $b(t) \uparrow \infty$ and a non-zero measure {$\nu(\cdot)\in \bM^*(\bC,\bO)$} such that as $t \to \infty$,}
\begin{equation}\label{eqn:ccc}
tP\left[ \frac{{\bZ }}{b(t)} \in \cdot \right] \stackrel{*}{\rightarrow} \nu(\cdot) \hskip 1 cm \text{in {$\bM^*(\bC, \bO)$.}}
\end{equation}
{When there is no danger of confusion,}  we sometimes use the
notation $\bM^*(\bO)$ to mean $\bM^*(\bC, \bO)$ {and sometimes
  abuse language and 
say the distribution is regularly varying on $\bO$.}}
\end{defn}

Definition \ref{def:closedconecomp} {implies}  there exists
$\alpha >0$ such that $b(\cdot)\in RV_{1/\alpha}$ and that  $\nu$
has the scaling {property: 
\begin{align}\label{eq:scaling} 
 \nu(c\;\cdot)=c^{-\alpha} \nu(\cdot), \hskip 1 cm c > 0.
\end{align}
 }This can be derived as in \cite[Theorem 3.1]{hult:lindskog:2006a}. 
We define standard multivariate regular variation, hidden
regular variation and the conditional extreme value model 
 in terms of Definition \ref{def:closedconecomp} and attempt to relate
 each to the way these ideas were first proposed
on modifications of compactified spaces.

Examples:
\begin{enumerate}
\item\label{eg:R+d} Let $\bC= [0,\infty)^d $ and $\bF=\{\bzero\}$ and
  $\bD=\bO=[0,\infty)^d \setminus \{\bzero\}$. Regular variation on
  $\bD$ is equivalent to regular variation defined in
  \eqref{eqn:existing_reg_var}  on $\E$.  The definition in
  \eqref{eqn:existing_reg_var}  precludes $\mu$ having mass on the
  lines through $\binfty$.
See Appendix \ref{subsec:equivRegVar}.
\item\label{eg:Dl} Let $d=2$, $\bC= [0,\infty)^2 $ and $\bF=\{(x,0):x\geq 0\}\cup
  \{(0,y), y\geq 0\} $ and $\bO=(0,\infty)^2$, the first quadrant with
  both the $x$ and $y$ axes removed.  This is the
  restriction to $[0,\infty)^2$ of the cone used in the definition of
  hidden regular variation in \citet{resnick:2002a}. For $d>2$, other
  examples of $\bF$ are in \citet{mitra:resnick:2010} and Subsection
  \ref{rem:relation_with_existing} provides a comparison
  between regular variation defined in \eqref{eqn:ccc} on
  $\bD^{(l)}=[0,\infty)^d\setminus \{\bx   \in
  [0,\infty)^d:x^{(l)}>0\} $ and regular variation defined using \eqref{eqn:existing_reg_var}
on $\E^{(l)}=[0,\infty]^d\setminus \{\bx
  \in \E: x^{(l)}>0\}$ where recall $x^{(l)}$ is the $l$th largest
  component of $\bx$. The two notions are equivalent provided there is
  no mass on $\E^{(l)} \setminus \bD^{(l)}$. 
  
\item 
Suppose $d=2$ and $\bC= [0,\infty)^d $ and $\bF=\{(x,0):x\geq 0\}.$
Then $\bO=\{(x,y) \in [0,\infty)^2: y>0\}=\bD_\sqcap$, the first
quadrant with the $x$-axis removed. This is the restriction to
$[0,\infty)^2$ of the cone used in the definition of the conditional
extreme value model \citep{heffernan:tawn:2004,
  heffernan:resnick:2007, das:resnick:2011}.
\item 
Suppose $d=2$ and $\bC= [0,\infty)^d $ and $\bF=\{(x,x):x\geq 0\}.$
Then $\bO$ is the first quadrant with the diagonal removed. This
example is suitable for discussing asymptotic full dependence
(\citet[page 294]{resnickbook:2008}, \citet[page 195]{resnickbook:2007})
and is
considered in Example
\ref{ex:asym_full_dep}.
\end{enumerate}

\subsection{Spectral measures, unit spheres and semi-parametric
  representations}\label{subsec:spectral}
Regular variation on $\E$ using the vague convergence definition as in
\eqref{eqn:existing_reg_var} allows a polar coordinate transformation
$\bx \mapsto (\|\bx\|, \bx/\|\bx\|) $. Assuming $b(\cdot) \in
RV_{1/\alpha}$, the limit measure has the scaling property and when
this is expressed in polar coordinates yields the version of
\eqref{eqn:existing_reg_var} 
$$tP[\frac{\|\bZ\|  }{b(t)} \in dx , \frac{\bZ}{\|\bZ\|} \in
d \ba \,]
\stackrel{v}{\to} \alpha x^{-\alpha -1}dx \,\times \, {S^*(d\ba)}
 $$
where $S^*$ is a finite measure on $\partial \aleph=\{\ba \in \E: \|\ba\|=1\}$,
the unit sphere. {Fixing  $S^*(\partial \aleph)=c$, we define
  $S(\cdot)=S^*(\cdot)/c$ which becomes a probability measure 
on $\partial\aleph$} {called the spectral or angular measure.} 
So in polar coordinates, the limit measure $\mu$ in 
\eqref{eqn:existing_reg_var}  has a semi-parametric {product} structure
{depending on the parameter $\alpha$ and the measure $S$.}

In $\E$, the unit sphere $\partial \aleph=\{\bx \in \E: \|\bx\|=1\}=\{\bx \in \E:
d(\bx, \bzero)=1\}$ is compact. {However, this may no longer be true
when moving to other sub-cones.}
 For instance in $(0,\infty]^2 $ the usual unit sphere is not
relatively compact. While the {polar coordinate transformation
  still allows this semi-parametric representation for other cones,} the analogue of $S$ is no
longer necessarily finite and this is a problem for inference.
{We explain next how to use a change of coordinates different from the
polar coordinate transformation which always produces a finite spectral
measure analogue.}  {\citet{heffernan:resnick:2007} and 
\citet{mitra:resnick:2010} consider alternatives to the
polar coordinate transformation that
 twist limit measures into a semi-parametric
form.}

{Proceed using the} context of  Definition
\ref{def:closedconecomp}. 
Assume $\bC, \bF$ and $\bO$ {are} defined as in Definition \ref{def:closedconecomp}
and define
 $\aleph_{\mathbb{O}} =  \{ \bx \in \mathbb{O}:
d(\bx, \mathbb{F}) \ge 1\}$ and
 $\aleph_{\mathbb{O}}$ is a subset of $\mathbb{O}$ bounded away
from $\mathbb{F}$ and $\cup_{\{ \theta > 0\}} \theta
\aleph_{\mathbb{O}} = \mathbb{O}$. From Definition
\ref{def:closedconecomp},  we have $0 < \nu(\aleph_{\mathbb{O}}) <
\infty$. The  scaling property \eqref{eq:scaling} implies
the scaling function $b(t)$ in \eqref{eqn:ccc}
can be chosen so that  so that
 $\nu(\aleph_{\mathbb{O}}) = 1$. Define the related set 
 $\partial \aleph_{\mathbb{O}} =  \{ \bx \in \mathbb{O}:
d(\bx, \mathbb{F}) = 1\}$.

Examples: (i) $\bC= [0,\infty)^d $ and $\bF=\{\bzero\}$ and
$d(\bx,\bF)=\|\bx\|$ and $\partial \aleph_\bO =\{\bx : \|\bx\|= 1\}.$
(ii) $d=2$, $\|\bx\|=x^1\vee x^2$, $\bC= [0,\infty)^2 $ and $\bF=\{(x,0):x\geq 0\}\cup
  \{(0,y), y\geq 0\} $ and $\bO=(0,\infty)^2$.  Then $\partial
  \aleph_\bO =\{\bx: d(\bx,\bF)=1\}=\{\bx: x^1 \wedge x^2=1\}.$
(iii) $d=2$, $\|\bx\|=x^1\vee x^2$,  and $\bC= [0,\infty)^d $ and $\bF=\{(x,0):x\geq 0\}.$
Then  $\partial
  \aleph_\bO =\{\bx: d(\bx,\bF)=1\} =\{(x,y): x\geq 0, y=1\}.$

{We} transform to an appropriate coordinate system {in which} 
the limit measure $\nu$ in \eqref{eqn:ccc} {is}
 a product of two components: a
one-dimensional Pareto measure and a 
probability measure defined on $\partial \aleph_{\mathbb{O}}$ called
the hidden spectral or angular 
measure.  To do this note
 two  properties of the distance function $d(\cdot, \mathbb{F})$: 
\begin{enumerate}[(i)]
\item Since $\mathbb{F}$ is a closed subset of $\mathbb{C}$,
  $d(\bx, \mathbb{F}) > 0$ for all $\bx \in \bO$ (else $\bx
  \in \bF$); and 
\item Since $\mathbb{F}$ is a cone, $\theta \cdot \mathbb{F} = \mathbb{F}$ for $\theta > 0$. Hence, $d(\theta \bx, \mathbb{F}) = d(\theta\bx, \theta \mathbb{F}) = \theta d(\bx, \mathbb{F})$, that is $d(\cdot, \mathbb{F})$ is homogeneous of order $1$.

\end{enumerate}

{A lemma is}
necessary for  the decomposition of limit measure $\nu(\cdot)$. For a
set $A \in [0,\infty)^d$ we set $(A)_1=\{x^1: \bx \in A\}.$

\begin{lem}\label{lem:smalldecomp}
Suppose  $h: \bO \mapsto (0, \infty) \times \partial \aleph_{\bO}$ is
a continuous bijection satisfying
\begin{enumerate}[(i)]
\item For every measurable $A \subset (0, \infty) \times \partial
  \aleph_{\bO}$ with $\overline {(A)_1} \cap \{ 0 \} = \emptyset$, 
 $h^{-1}(A)$ is
  bounded away from $\mathbb{F}$. 
\item For every measurable $B \subset \bO$ with $B$ bounded away from
  $\mathbb{F}$,  $\overline{{(h(B))_1}} \cap \{ 0 \}
  = \emptyset$. 
\end{enumerate}
Then, the following statements are equivalent:
\begin{enumerate}[(i)]
\item  As $t \to \infty$,
\begin{equation*}
\mu_t(\cdot) \stackrel{*}{\rightarrow} \mu(\cdot) \hskip 1 cm \text{in } \bM^*(\bC, \bO).
\end{equation*}

\item For all measurable $A  \subset (0, \infty) \times \partial
  \aleph_{\bO}$ such that $\overline{(A)_1} \cap \{ 0 \} = \emptyset$ and
  $\mu \circ h^{-1} (\partial A) = 0$, 
\begin{equation*}
\mu_t \circ h^{-1}(A) \rightarrow \mu \circ h^{-1}(A).
\end{equation*}
\end{enumerate}

\end{lem}
\begin{proof} The proof  follows the steps of the proof of Theorem 2.5
  of \citet{hult:lindskog:2006a}.
\end{proof}
Now by applying 
Lemma \ref{lem:smalldecomp} with
$h : \bx \mapsto \left(d(\bx,
  \mathbb{F}), {\bx}/{d(\bx, \mathbb{F})} \right)$, we are able to decompose $\nu$ as follows.

\begin{prop}\label{prop:smallconedecomp}
Regular variation on $\bO$ as given in \eqref{eqn:ccc} is equivalent to 
\begin{equation}\label{eqn:smallcone_decomp}
tP\left[ \left( \frac{d({\bZ }, \mathbb{F})}{b(t)},
    \frac{{\bZ }}{d({\bZ }, \mathbb{F})} \right) \in A \right]
\rightarrow {c}\nu_{(\alpha )}   \times S_{\bO}(A) 
\end{equation} 
for all measurable $A \subset (0, \infty) \times \partial
\aleph_{\bO}$ such that $\overline{(A)_1} \cap \{ 0 \} = \emptyset$ and
$\nu\circ h^{-1} (\partial A) = 0$ where $c>0$, 
 $S_{\bO}(\cdot)$ {is} a probability measure on $\partial \aleph_{\bO}$
and $\nu_{(\alpha )}(\cdot)$ is the Pareto measure given by
$\nu_{( \alpha )}((x, \infty) ) = x^{-\alpha}$ for $x >0$. Call
 $S_{\bO}(\cdot)$  the  spectral measure on
$\bO$; it  is related to $\nu(\cdot)$ by the relation 
\begin{equation}\label{eqn:s1andnu1}
S_{\bO}(\Lambda) = \frac{\nu \left( \left\{ \bx \in \bO: d(\bx, \mathbb{F}) \ge 1, \frac{\bx}{d(\bx, \mathbb{F})} \in \Lambda \right\} \right)}{{\nu\left(  \bx \in \bO: d(\bx, \mathbb{F}) \ge 1\right)}}.
\end{equation}
\end{prop}

Since $(0,\infty)\times \partial \aleph_{\bO}$ is not a cone, we have
not phrased the convergence in \eqref{eqn:smallcone_decomp} as $\bM^*$
convergence as in \eqref{eqn:ccc}. To do so would require {more}
reworking of
the convergence theory in \citet{hult:lindskog:2006a}.}

\begin{cor}\label{cor:decomp2statements}
 The convergence in \eqref{eqn:ccc} is equivalent to the following two conditions:
\begin{enumerate}[(i)]
\item The distribution of $d({\bZ }, \mathbb{F})$ is regularly
  varying on $(0, \infty)$ following Definition
  \ref{def:closedconecomp} with $\bC = [0, \infty)$. 

\item The conditional distribution of ${\bZ }/d({\bZ },
  \mathbb{F})$ given $d({\bZ }, \mathbb{F}) > t$, converges weakly,
\begin{equation*}
P\left[ {\bZ }/d({\bZ }, \mathbb{F}) \in \cdot \hskip 0.1 cm
  \big{|} \hskip 0.1 cm  d({\bZ }, \mathbb{F}) > t \right]
\Rightarrow S_{\bO}(\cdot) \qquad (t\to\infty). 
\end{equation*}
\end{enumerate}
\end{cor}

\begin{rem}\label{rem:smalldecomp}
{\rm
We make a few remarks about Proposition \ref{prop:smallconedecomp}.
\begin{enumerate}[(i)]
\item \label{rem:why_not_compact} On the role of the distance function: Proposition
  \ref{prop:smallconedecomp} {emphasizes} that the 
  spectral probability measure $S_{\bO}(\cdot)$ is dependent on the
  choice of distance function
  $d(\cdot,\cdot)$.
Corollary \ref{cor:decomp2statements} allows us to use the
distance function   $d(\cdot,\cdot)$
 to detect 
regular variation on $\bO$; see Section \ref{sec:detection_estimation}. 
However, extending the distance function 
to a compactified space such as
  $[0, \infty]^d$ is difficult and  this provides another  reason why
  we deviated from the standard {discussion of regular variation using} compactified spaces and
  vague convergence.

\item Connections to prior treatments:
\begin{enumerate}
\item Proposition 3.1 of \cite{mitra:resnick:2010} decomposes the
  limit measure {$\mu^{(l)}(\cdot)$ on $\E^{(l)}$ (see
    \eqref{eqn:existing_hidden_defn} below)  by applying the transformation
  $T : \bx \mapsto \left( x^{(l)}, {\bx}/{x^{(l)}}
  \right),$ where $x^{(l)}$ is the $l$-th largest component of
  $\bx$. If we choose $\mathbb{C} = [0, \infty)^d$, 
$\mathbb{F} = \{
  \bx \in \mathbb{C} : x^{(l)} = 0 \}$ and
$\bO = 
\{   \bx \in \mathbb{C} : x^{(l)} > 0 \}$, and choose the
  $L_{\infty}$-norm when defining
$d(\cdot, \cdot)$, then $d(\bx, \mathbb{F}) = x^{(l)}$ and 
{our}
Proposition \ref{prop:smallconedecomp} gives a  version 
 of Proposition 3.1 in \citet{mitra:resnick:2010}. 

\item 
For considering 
  regular variation on the cone ${\E_\sqcap :}=[0, \infty] \times (0, \infty]$,
  \citet[Proposition 4]{heffernan:resnick:2007}  give a 
  decomposition of {their} limit measure $\mu_*(\cdot)$ by applying the
  transformation $T: (x,y) \mapsto (y, x/y)$,  If
  we choose $\mathbb{C} = [0, \infty)^2$,
$\mathbb{F} = [0, \infty) \times \{ 0 \}$ and 
 $\bO = [0, \infty) \times
  (0, \infty)=:{\D_\sqcap}$,  and define
  $d(\cdot, \cdot)$ using the
  $L_{\infty}$-norm, 
then $d\bigl((x,y), \mathbb{F}\bigr) = y$ and our Proposition
\ref{prop:smallconedecomp} 
connects with   \citet[Proposition 4]{heffernan:resnick:2007}.

\item Proposition \ref{prop:smallconedecomp}
also relates to the usual polar coordinate characterization of
multivariate regular variation on $\E$ as in \citet[page
173]{resnickbook:2007}. Set $\bF=\{ 
{\bzero } \}$ which is a closed cone {in $\bC=[0,\infty)^d$ and 
$d(\bx, \{ {\bzero } \}) = ||\bx||$. See also Proposition
\ref{prop:bigconedecomp}.
}}
\end{enumerate}
\end{enumerate}
}
\end{rem}

\subsection{Hidden regular variation}
As in Definition \ref{def:closedconecomp}, consider 
$\bF \subset \bC \subset [0,\infty)^d$ with $\bF$ and $\bC$ closed
cones containing $\{\bzero\}$. Suppose $\bF_1$ is another subset of
$\bC$ that is a closed cone containing $\{\bzero\}$.  Then 
$\bF \cup \bF_1$ is also a closed cone containing $\{\bzero\}$.
Set $\bO= \bC\setminus \bF$ and
$\bO_1= \bC \setminus (\bF \cup \bF_1)$.

 \begin{defn}\label{hidden_defn1}
 {\rm{
 The distribution of  a random vector ${\bZ } \in \mathbb{C}$ that
 is regularly varying on $\mathbb{O}$ with scaling function $b(t)$ in
 \eqref{eqn:ccc} 
possesses  hidden regular variation (HRV)
 on $\mathbb{O}_1$ if 
\begin{enumerate}
\item The distribution of $\bZ$ is also regularly varying on $\bO_1$
  with scaling function $b_1(t) $ and limit measure $\nu_1$ and 
\item $b(t)/b_1(t) \to \infty$ as $t\to\infty$.
\end{enumerate}
}}
 \end{defn}
 
Observe that the condition $b(t)/b_1(t)\to \infty$ implies $\nu$ puts
zero mass on $\bO_1$. (See \citet{resnick:2002a}.)
 From Definition \ref{hidden_defn1} {and \eqref{eq:scaling}}, it
 follows that there exists $\alpha_1 \ge \alpha$ such that $b_1(\cdot)
 \in RV_{1/\alpha_1}$ and on $\mathbb{O}_1$, the limit measure
 $\nu_1(\cdot)$ in \eqref{eqn:ccc} satisfies
 the scaling property 
\begin{equation}\label{scaling_c1}
\nu_1 (c \cdot) = c^{-\alpha_1}\nu_1(\cdot), \hskip 1 cm c > 0.
\end{equation}

{Example: Let $\bZ =(Z^1,Z^2)$ be iid unit Pareto random variables.
Then the distribution of $\bZ$ is regularly varying on $\bD$ with
$b(t)=t$ and possesses HRV on $(0,\infty)^2$ with $b_1(t)=\sqrt t.$
Somewhat more generally, if ${\bZ } = (Z^1, Z^2)$ has regular
variation on $\bD= [0, \infty)^2 \setminus \{{\bzero }\}$,
then in  the presence of asymptotic independence, HRV offers non-zero 
estimates of joint tail probabilities $P[ Z^1 > x, Z^2 >
y]$ for large thresholds $x, y > 0$, whereas regular variation on
$\bD$ estimates joint tail probabilities as zero.
Other examples are considered in Section \ref{sec:examples}.}

\begin{rem}\label{remark: hiddendefn1}
{\rm{
We make a few remarks on Definition \ref{hidden_defn1}. 

\begin{enumerate}[(i)]
\item \label{rem:anysubcone1} There is flexibility in
    choosing $\bO_1$ and this flexibility is useful for 
  defining HRV on a   sequence of cones. 
Cones can be  chosen based on the risk regions whose probabilities are required.
For example, if $d = 2$, we choose the
    cones $\bD$ and $\bD^{(2)} := \{\bz  \in [0, \infty)^d:
 z^1\wedge   z^2 > 0 \}$ if we need the probability that components of
 the risk vector
simultaneously exceed thresholds.

\item 
{\it{{Differences} with
      existing notions of hidden regular variation}}: 
Previous considerations of HRV relied on vague convergence and
compactification and were applied to specific choices of cones.
\citet{resnick:2002a}, \citet{heffernan:resnick:2005}, and \citet{maulik:resnick:2005} consider HRV
  on $(0,\infty]^2$ and \citet{mitra:resnick:2010} consider the cone
  $\E^{(l)} = \{ \bx \in [0, \infty]^d: x^{(l)} > 0 \}$. 
The choice of these specific cones {may} not provide sufficient
flexibility and generality. For example,  to estimate $P[|X-Y|>x]$ when asyptotic full dependence is
present, such cones considered previously are of no help. See Example \ref{ex:asym_full_dep}.   
\end{enumerate}
}}
\end{rem}

{\subsubsection{Where to seek HRV}\label{subsubsec:where?}
Suppose the distribution of $\bZ$ is regularly varying on $\bO$ and
that the limit measure $\nu$ in \eqref{eqn:ccc} gives zero mass to a
subset $R$ of $\bO$. Using the asymptotic property of regular
variation to estimate $P[\bZ \in t R]$ {for large $t$,} means such an estimate is $0$.
So we seek another regular variation on a subset of $R$ which is of
lower order.

Thus, when seeking HRV our focus is on subsets of $\mathbb{O}$ where the
limit measure $\nu(\cdot)$ gives zero mass.
 A systematic way to find HRV is facilitated by the
following simple remark.}

\begin{prop}\label{prop:support}
In Definition \ref{def:closedconecomp}, the support of the limit
measure $\nu$ is a closed cone $  F_{\nu} \subset \bC$ containing
$\bzero$.
\end{prop}

\begin{proof}
Let $\supp (\nu)$ denote the support of $\nu$.  By definition, $\supp (\nu)$
is closed. Let $\bx \in \supp (\nu)$ and we show for $t>0$ that $t\bx
\in \supp (\nu)$. For small $\delta $, by \eqref{eq:scaling}
\begin{align*}
\nu \bigl( (t\bx-\delta \bone, &t\bx +\delta \bone )  \cap \bC \bigr)
=\nu \Bigl( t\bigl( (\bx-\frac{\delta}{t} \bone, \bx +\frac{\delta }{t}
\bone )  \cap \bC \bigr) \Bigr)\\
=&t^{-\alpha}  \nu \bigl( (\bx-\frac{\delta}{t} \bone, \bx +\frac{\delta }{t}
\bone )  \cap \bC \bigr)>0
\end{align*} since $\bx \in \supp (\nu)$.
\end{proof}

So $\bF_\nu=\supp (\nu)$ is a union of rays emanating from the origin. It is
also true that 
$$\supp (\nu) = \{t \cdot \supp (S_\bO ),\, t\geq 0\}.$$

When seeking HRV on a cone smaller than $\bO$, we conclude that
$$\bO_\nu := \bC \setminus (\bF \cup  \bF_\nu)$$
is the  largest possible sub-cone of
$\bO$ where we might find a different regular variation.
In practice, guided by the type of risk region whose probability we
need to estimate, we find a closed cone $\bF_1 \subset \bC$ containing
$\bzero$ such that
$\bF_1 \supset \bF_\nu$ and set $\bO_1=\bC \setminus (\bF \cup \bF_1)$ and
then seek regular variation on $\bO_1$. Possibly, but not necessarily
$\bF_1= \supp (\nu)$. Examples are in \citet{mitra:resnick:2010} and
Section \ref{sec:examples}.

\subsubsection{Regular variation on a sequence of cones}\label{subsubsec:sequence}
Having found regular variation on $\mathbb{O}$ with HRV on
$\mathbb{O}_1$, we ask: should we stop here? There might be a subcone
$\mathbb{O}_2$ of $\mathbb{O}_1$, where $\nu_1(\cdot)$ gives zero mass
and hence, there might exist a different regular variation on
$\mathbb{O}_2$.  

To proceed further, as before remove the support of $\nu_1$ from $\bO_1$ and consider the set
$$ \bO_{\nu_1} := \bO_1 \setminus \supp ( \nu_1 )$$
and  $\bO_{\nu_1}$ is the largest possible sub-cone of $\bO_1$ where
we might find a different regular variation. So we choose
$\bF_2 \supset \supp (\nu_1) $ and set $\bO_2 = \bC \setminus (\bF \cup
\bF_1 \cup \bF_2 )$ and seek regular variation on $\bO_2$ with scaling
function $b_2(t)$ such that $b_1(t)/b_2(t) \to \infty$ as $t
\to\infty$. This last condition guarantees the regular variation on
$\bO_2$ is of lower order than the regular variation on either $\bO$
or $\bO_1 $ and hidden from both higher order regular variations. This
process of discovery is continued as long as on each new cone regular
variation is found. Example \ref{ex:infinite_sequence} shows this
discovery process may lead to an infinite sequence of cones. 

From our definition of HRV, at each stage {of the discovery
  process}
we have some flexibility in
choosing the next cone where we seek HRV. 
Example \ref{ex:infinite_sequence} shows it may be impractical to 
analyze HRV on every possible cone as the discovery process may lead
to an infinite
sequence of cones.  A more
practical approach is to decide on a particular finite sequence
of cones based on the risk regions {of  interest};
 see Remark
\ref{remark: hiddendefn1}{(\ref{rem:anysubcone1})}. 
For example, if we
are interested in estimating joint tail probabilities, we might consider only
the sequence of cones $\mathbb{D} =  [0, \infty)^d \setminus \{
{\bzero } \} \supset \mathbb{D}^{(2)} = \{ \bx \in [0, \infty)^d:
x^{(2)} > 0\} \supset \cdots \supset \mathbb{D}^{(d)} = \{ \bx
\in [0, \infty)^d: x^{(d)} > 0\}$; cf. \citet{mitra:resnick:2010}.

\section{{Remarks on} other models of multivariate regular
  variation}\label{sec:CEVandNon-stanMRV_connection} 
Despite the fact that most common examples in heavy tail analysis
start by analyzing convergence on the cone $[0,\infty)^d \setminus
\{\bzero\}$, this need not always be the case. 
For example, the standard case of the conditional extreme value (CEV) model
\citep{heffernan:resnick:2007, das:resnick:2011}, is regular variation
\eqref{eqn:ccc} with $b(t)=t$ 
on the cone $\bD_{\sqcap}:=\bO=[0, \infty) \times (0, \infty)$ with $\bC = [0,
\infty)^2$ and $\bF=[0,\infty) \times \{0\}$.

\subsection{The CEV model}\label{subsec:cev}

The conditional extreme value (CEV) model, suggested in 
\citet{heffernan:tawn:2004}, is an alternative to  classical
multivariate extreme value theory (MEVT). In contrast with  
classical  MEVT which implies all marginals are in a maximal domain of
attraction,
in the CEV model only a particular subset of the random vector is assumed
to be in a maximal domain of attraction.
For convenience, restrict attention to $d=2$.


The CEV model as formulated in 
\citet{heffernan:resnick:2007, das:resnick:2011} allows variables to
  be centered as well as scaled. To make comparison with models of
  regular variation on the first quadrant easy, we recall the vague
  convergence definition using only scaling functions. 

\begin{defn}\label{CEV_defn}
 {\rm{
 Suppose $\bZ:=(\xi,\eta) \in \R_{+}^2$ is a random vector and there exist functions $a_1(t),a_2(t)>0$ and a non-null Radon
 measure $\mu$ on Borel subsets of $\bE_{\sqcap}:=[0,\infty]\times(0,\infty]$ such that in $\bM_{+}(\bE_{\sqcap})$
\begin{align}\label{eqn:CEV}
 & t P\left[\left(\frac{\xi}{{a_1}(t)},\frac{\eta}{{a_2}(t)}\right) \in \; \cdot\;\right] \cnvg \mu(\cdot).
\end{align}
Additionally assume that $\mu$ satisfies the following non-degeneracy conditions:
\begin{itemize}
 \item[(a)] $\mu([0,x]\times(y,\infty])$ is a non-degenerate distribution in $x$,
 \item[(b)] $\mu([0,x]\times(y,\infty]) < \infty$.
\end{itemize}
Also assume that 
\begin{itemize}
 \item[(c)] $H(x):=\mu([0,x]\times(1,\infty])$ is a probability distribution. 
\end{itemize}
In such a case
 $\bZ$ satisfies a conditional extreme value model and we write $\bZ \in CEV({a_1,a_2})$.
}}
 \end{defn}

The general CEV model, provided the limit measure is not a product,
 can be standardized to have standard regular variation on
the cone $\E_{\sqcap}$  \citep[pg. 236]{das:resnick:2011}.
Following the theme of Examples \ref{eg:R+d} and \ref{eg:Dl} at the
beginning of Section \ref{sec:rvandhrv}, if no
mass exists on the lines through $\binfty$, and {$a_1=a_2$}, then the
vague convergence definition of the CEV model on $\E_\sqcap $ is the same as the
$\bM^*$ definition on regular variation $\bD_\sqcap$.

The issue of mass on the lines through $\binfty$ is significant since
if mass on these lines is allowed, there is a statistical
identifiability problem in the sense that  in $\E_\sqcap$ it is possible to have two
different limits in \eqref{eqn:CEV} under two different normalizations
and under one normalization, there is mass on the lines through
$\binfty$ and with the other, such mass is absent.
 Restricting to $\bD_\sqcap$ resolves the {identifiability} problem
as in this space, limits are unique. 
See Example \ref{ex:CEV_id}.

\subsection{Non-standard regular variation}\label{subsec:nonstdrv}
Standard multivariate regular variation on $\E$ requires the same
normalizing function to scale all  components and is a convenient
starting place for theory but  unrealistic for applications
as it makes all one dimensional marginal distributions tail equivalent.
Non-standard regular variation \citep[Section 6.5.6]{resnickbook:2007}
allows different normalizing functions
for vector components and hence permits
each marginal distribution tail to have a different tail index. When
the components of the risk vector have different tail indices,
non-standard regular variation is sensitive to the different tail
strengths. On $\E$ or $\bD$, non-standard regular variation takes the form
\begin{equation}\label{eqn:cccc}
tP[ (Z^i/{a_i}(t), i=1,\dots,d) \in \cdot ] \to \nu(\cdot) ,
\end{equation}
for scaling functions {$a_i(t) \uparrow \infty$}, $ i=1,\dots,d$, where
convergence is vague for $\E$ and in $\bM^*$ for $\bD$.
If convergence is in $\E$ and there is no mass on the lines through
$\binfty$, the difference between convergence in $\E$ and $\D$ evaporates.

In cases where {$a_i(t)/a_{i+1} (t) \to 0$}, it is sometimes possible to compare the
information in non-standard regular variation with what can be
obtained from HRV. Sometimes HRV provides more detailed information. Consider the following
{
\begin{ex} \label{eg:HRVvsNonST}
{\rm{
Suppose $X_1, X_2, X_3 $ are independent random variables where $X_1$ is Pareto($1$), $X_2$ is Pareto($3$) and $X_3$ is Pareto($4$). Further assume
that $B$ is a Bernoulli$(1/2)$ random variable independent of  $X_1, X_2, X_3 $. Define
\[ \bZ=(Z^1,Z^2):= B(X_1,X_3) + (1-B)(X_2,X_2).\]

Non-standard regular variation on $\E$ or $\D$ is given by

$$tP\left[ \left(\frac{Z^1}{t}, \frac{Z^2}{t^{1/3}} \right) \in dx \,dy\right]
\to \frac 12 x^{-2}dx \cdot \epsilon_0 (dy)+ \frac 12\epsilon_0(dx)\cdot 3y^{-4}dy,$$
where $\epsilon_0$ indicates the point measure at $0$ and the limit
measure concentrates on the two axes. Now  HRV can also be sought under such non-standard regular variation and  we get 
$$ tP\left[ \left(\frac{Z^1}{t^{3/7}}, \frac{(Z^2)^3}{t^{3/7}} \right) \in dx \,dy\right] {\to} \frac 12 x^{-2} dx\, 4y^{-5} dy$$
on the space $\E^{(2)}$. Note that, this form of regular variation completely ignores the existence of the completely dependent component
of $\bZ$ given by $(X_2,X_2)$. Alternatively, if we pursue regular variation and HRV on a sequence of cones as defined in this paper then
 we observe the following convergences as $t \to \infty$:
\begin{enumerate}
 \item On $\D$, we have
\begin{align*}
tP\left[ \left(\frac{Z^1}{t},\frac{Z^2}{t} \right) \in dx \,dy\right] \to \frac1{2}x^{-2}dx\, \epsilon_0(dy).
\end{align*}
\item In the next step on $\D \setminus {\{\text{x-axis} \}}$, we have
\begin{align*}
tP\left[ \left(\frac{Z^1}{t^{1/3}},\frac{Z^2}{t^{1/3}} \right) \in dx \,dy\right] \to \frac1{2}3x^{-4}dx\, \epsilon_x(dy). 
\end{align*}
\item Next, on $\D \setminus \left[{\{\text{x-axis} \}} \cup \{ \text{diagonal}\}\right]$, we have
\begin{align*}
tP\left[ \left(\frac{Z^1}{t^{1/4}},\frac{Z^2}{t^{1/4}} \right) \in dx \,dy\right] \to \frac1{2}\epsilon_{0}(dx)\,4y^{-5}dy.  
\end{align*}
\item Finally, on $\D \setminus \left[{\{\text{x-axis} \}} \cup \{ \text{diagonal}\} \cup {\{\text{y-axis} \}} \right]$, we have
 \begin{align*}
tP\left[ \left(\frac{Z^1}{t^{1/5}},\frac{Z^2}{t^{1/5}} \right) \in dx \,dy\right] \to \frac1{2}x^{-2}dx\,4y^{-5}dy.  
\end{align*}
Clearly, this analysis captures the structure of ${\bf{Z}}$ better than what non-standard regular variation along with HRV in the classical set-up.

\end{enumerate}

}}\end{ex}}
\COM{
\begin{ex} \label{eg:HRVvsNonST}
{\rm{
Suppose $\xi \independent \eta$ and $\xi$ is Pareto(1) and $\eta$ is
Pareto(2). Non-standard regular variation on $\E$ or $\bD$ is
$$tP[ \bigl(\frac{\xi}{t}, \frac{\eta}{\sqrt t} \bigr) \in dx \,dy]
\to x^{-2}dx \cdot \epsilon_0 (dy)+ \epsilon_0(dx)\cdot 2y^{-3}dy,$$
where $\epsilon_0$ indicates the point measure at $0$ and the limit
measure concentrates on the two axes. Comparing with 
\eqref{eqn:cccc}, in this example
$b_1(t)=t$ and $b_2(t)=\sqrt t$.
Standard regular variation coupled with HRV reveals the following
information in stages:
\begin{enumerate}
\item Using the heavier normalization $b_1(t)=t$,
$$ tP[\bZ/t \in dx\, dy] \to x^{-2}dx \cdot \epsilon_0 (dy),$$
on $\E$ or $\bD$ since $b_1(t)$ is too heavy for the $\eta $
variable. This limit measure concentrates on the horizontal axis.
\item Peeling away the horizontal axis, using normalization
  $b_2(t)=\sqrt t$ we get (remember $\eta$ is Pareto(2)),
$$tP[ \bZ/{\sqrt t} \in dx\, dy] \to \epsilon_0 (dx) \cdot 2y^{-3}dy,$$ 
on $\E_\sqcap$ or $\bD_\sqcap$ and the limit concentrates on the
vertical axis of $\bD_\sqcap$.
\item HRV allows us to go further. Continue peeling away parts of the state
  space and focussing where the limit measure in the previous step
  put zero mass. So for this example we peel away the vertical axis and use
  $b_3(t)=t^{1/3}$ to get
$$tP[\bZ/t^{1/3} \in dx\,dy] \to x^{-2}dx \cdot 2y^{-3}dy,$$
on $(0,\infty]^2$ or $(0,\infty)^2$.
\end{enumerate}
}}\end{ex}
}
Thus, in Example \ref{eg:HRVvsNonST}, HRV provides more information than non-standard
regular variation. {T}his is  not always the case {and}
sometimes non-standard regular variation is better suited to explain{ing}
the structure {of} a model.

\begin{ex}\rm{\citep[Example 6.3]{resnickbook:2007}}\label{eg:HRV_blows}
{\rm{
Suppose $X$ is a standard Pareto
  random variable and   $\bZ = (X,X^2)$.
Using the obvious non-standard scaling for  the coordinates leads to 
\begin{align}
 tP\left[\left(\frac{X}{t},\frac{X^2}{t^2}\right)  \in  \cdot
   \,\right] \cnvg \nu_{(1)}\circ T^{-1}, ~~~~(\text{as } t \to
 \infty) \label{eqn:explain}
\end{align}
where $T:(0,\infty]\to (0,\infty]\times(0,\infty]$ is defined by
$T(x)=(x,x^2)$ and $\nu_{(1)} (dx)=x^{-2}dx,\,x>0.$ The limit measure
concentrates on $\{(x,x^2): x>0\}$. Using 
 the same  normalization on both coordinates is not so revealing. With
 the heavier normalization, we have on {$\bM^*([0, \infty)^2, [0,\infty)^2\setminus\{\bzero\}$)},
\begin{align}
 tP[\bZ/{t^2} \in \cdot
   ] \stackrel{*}{\to} \epsilon_0\times \nu_{(1)}, ~~~~(\text{as } t \to
\, \infty) \label{eqn:noexplain1} 
\end{align}
{and} with the lighter normalization and $x>0,\,y>0$,
$$P[X/t>x, X^2/t >y] \to x^{-1},\qquad \forall y>0,$$
so that 
\begin{align}
 tP[\bZ/{t}  \in \cdot] \cnvg
  \nu_{(1)}\times \epsilon_\infty,\quad \text{ on}\,\E. \label{eqn:noexplain2}
\end{align}
Neither \eqref{eqn:noexplain1} nor \eqref{eqn:noexplain2} approach the
explanatory power of \eqref{eqn:explain}. {Moreover, since our modeling excludes points at infinity, the convergence in \eqref{eqn:noexplain2} is not equivalent to any $\bM^*$-convergence.}

}}\end{ex}

\section{Examples}\label{sec:examples}
The definition of HRV given in this paper {changes somewhat the definition of
convergence but more importantly allows general cones compared} with the
existing notion of HRV \citep{resnick:2002a, mitra:resnick:2010}; see
Remark \ref{remark: hiddendefn1} and equation \eqref{rem:relation_with_existing} of
this paper. Here we provide examples to illustrate  use and highlight
subtleties.
The examples give risk sets, whose probabilities can
 be approximated  by using our general concept of HRV, and not the existing notion.
See also   \citep{mitra:resnick:2010}.

\begin{ex} \label{example:div_risk}
{\rm{ {\bf{Diversification of risk:}} 
{Suppose, we invest in two
    financial instruments $I_1$ and $I_2$ and  for a given time
    horizon future losses
    associated with each unit of $I_1$ and $I_2$ are $\xi$ and $\eta$
    respectively. }
Let ${\bZ } = \left(\xi, \eta \right)$ be regularly
    varying on $\bO = \mathbb{D} = [0, \infty)^2 \backslash \{
    {\bzero } \}$ with limit measure $\nu(\cdot)$. 
 {W}e earn risk premia of $\$l_1$ and
    $\$l_2$ for investing in each unit of $I_1$ and $I_2$. 

Suppose  we invest in 
$a_1$ units of $I_1$ and $a_2$ units of $I_2$. 
A possible risk
measure for this portfolio is 
$P[a_1\xi >x, a_2\eta >y]$ for two large thresholds $x$ and $y$ and this
risk measure  quantifies tail-dependence of $\xi$ and $\eta$. 
{The greater the tail-dependence between $\xi$ and $\eta$, the greater
should be our
reserve capital requirement so that we guard against the extreme}
situation that 
both investments go awry. For this risk measure, the best
circumstance is  if risk contagion is absent; that is,  $\xi$ and $\eta$ are asymptotically
independent {so}  that the measure $\nu(\cdot)$ is concentrated on the
axes \cite[page 192]{resnickbook:2007} 
since then $P[a_1\xi > x, a_2\eta > y]$ is estimated to be zero if HRV is absent
 and even if HRV exists on the cone $(0,
\infty)^2$ according to Definition \ref{hidden_defn1}, {the
  estimate of}
$P[a_1\xi > x, a_2\eta > y]$ for  large thresholds $x$ and $y$
{should be} small compared to the case where $\xi$ and $\eta$ are
not asymptotically independent. 

In place of asymptotic independence, suppose instead 
that  ${\bZ }$ is regularly varying on
$\mathbb{D}$ as in Definition \ref{def:closedconecomp} with limit
measure $\nu(\cdot)$ and $\nu(\cdot)$ has support $\left\{(u,v)
  \in \mathbb{D} : 2u \le v \right\} \cup \left\{(u,v) \in
  \mathbb{D} :u \ge 2v \right\}$ {so that
 $\nu(\cdot)$ puts zero mass on}
\begin{align}\label{eqn:CONE}
\text{CONE} :&=\left\{(u, v) \in \mathbb{D} : 2u > v, 2v >
  u \right\}\nonumber\\
&=\bigcup_{x>0, y>0} \{(u,v): 2u-v>x,2v-u>y\}   .\end{align}
We can {still} build a portfolio
of two 
financial instruments that have asymptotically independent risks as follows.

Define two new financial instruments $W_1 = 2I_1 -
I_2$ (buy two units of $I_1$ and sell a unit of $I_2$, assuming such
transactions are allowed) and $W_2 = 2I_2 - I_1$. {The loss associated
with $W_1$ is $L_{W_1}:=2\xi -\eta$ and the loss for $W_2$ is $L_{W_2}:=2\eta -\xi$.}
We earn the same risk premia {$a_1l_1+a_2 l_2$}  in the following two cases:
\begin{enumerate}[(i)]
\item \label{old_investment}  invest in $a_1$ units of $I_1$ and $a_2$ units of $I_2$, {or}
\item \label{new_investment}  invest in $c_1 = (2a_1 + a_2)/3$ units of $W_1$ and $c_2 = (a_1 + 2a_2)/3$ units of $W_2$. 
\end{enumerate}

A measure of risk of the portfolio in (i) is $P[a_1\xi >x, a_2 \eta
>y]$ and since as $t \to \infty$,
$$tP[a_1 \xi/b(t) >x,  a_2 \eta /b(t) >y] \to \nu\bigl( ((x,y),
\binfty))>0,$$
since asymptotic independence is absent, we expect the risk probability
to be not too small.  However, the risk measure for (ii) should be
rather small as we now explain. Let $T(u,v)=(2u-v,2v-u)$ and the risk
measure for (ii) is
$$P[c_1 (2\xi-\eta)>x, c_2
(2\eta-\xi)>y]
=P[T(\xi, \eta)> (x/c_1, y/c_2)].$$
Since as $t \to \infty$,
$$tP[T(\xi, \eta)> (b(t)x/c_1, b(t)y/c_2)]\to \nu\Bigl\{ T^{-1}\Bigl(
\bigl((x/c_1,y/c_2),\binfty \bigr) \Bigr)               \Bigr\}$$
and
$$ T^{-1}\Bigl(
\bigl((x/c_1,y/c_2),\binfty \bigr) \Bigr) \subset \text{CONE},$$
the risk measure for (ii) is small and the losses are asymptotically
independent.
 Hence, investment portfolio
\eqref{new_investment} above achieves more diversification of risk
than portfolio \eqref{old_investment}, and both 
earn the same amount of risk premium. 
 
How do we  construct a risk vector $\bZ$ whose distribution is regularly
varying and whose limit measure $\nu$ concentrates on $\bD\setminus
\text{CONE}$? 
Suppose $R_1,R_2, U_1,U_2,B$ are independent and $R_1$ is Pareto($1$),
$R_2$ is Pareto($2$), $U_1 \sim U\bigl( (0,1/3) \cup (2/3, 1)\bigr)$,
$U_2 \sim U(1/3,2/3)$ and $B $ is Bernoulli
with values $0,1$ with equal probability. Define
$$ \bZ=BR_1(U_1,1-U_1)+(1- B)R_2(U_2,1-U_2) .$$
Because of Proposition \ref{prop:smallconedecomp}  applied with the
$L_1$ metric, $\bZ$ is regularly varying on $\bD$ with index 1 and
angular measure concentrating on $[0,1/3]\cup [2/3,1]$. Hence the
limit measure $\nu $ concentrates on $\bD \setminus \text{CONE}$. At
scale $t^{1/2}$, $\bZ$
is also regularly varying on $\text{CONE}$ with uniform angular measure.

 }}
\end{ex}

\begin{ex}\label{ex:asym_full_dep}
{\rm{ 
{\bf Asymptotic full dependence:} 
{For convenience, restrict  {this example} to $d = 2$. }
HRV \citep{resnick:2002a} {was} designed to deal
with asymptotic independence \cite[page
322]{resnickbook:2007}
 where {$\nu(\cdot)$}  concentrates on the
axes. {For} asymptotic full dependence,
the limit measure {$\nu(\cdot)$} concentrates on the diagonal $
{\text{DIAG}:}=\{
(z, z) : z { \geq 0} \}$ \cite[page 195]{resnickbook:2007} {and}
previous treatments could not  deal with this or
{related} degeneracies where the}
limit measure ${\nu}(\cdot)$ concentrate{s} on a finite number of
rays other than the axes. 

{Consider} the following example. Suppose, 
 $X_1, X_2, X_3$ are
iid with common distribution Pareto($2$). Let $B$ be a Bernoulli
random variable independent of $\{ X_i: i = 1, 2, 3 \}$ and $P[B =0] =
P[B =1] = 1/2$. Construct the random vector ${\bZ }$ as 
$${ {\bZ } = {(\xi, \eta)} =  B({(X_1)}^2, {(X_1)}^2) +
  (1-B)(X_2, X_3)}$$
{and}
 ${\bZ }$ is regularly varying on $\bO = \mathbb{D}$ with scaling
 function $b(t) = t$ and limit measure $\nu$
{where,} 
$$
\nu \left( [{\bzero }, (u,v)]^c \right) = \frac1{2} {\left(u \wedge
    v\right)}^{-1}
\qquad  (u,v) \in \mathbb{D}.$$
The measure $\nu(\cdot)$ {concentrates on $\text{DIAG}$} and satisfies 
$\nu (\{ (u, v) \in \mathbb{D}: |u-v| > x
\} ) = 0$,  for $x > 0$ and 
as a result, {we estimate} {as zero} risk probabilities like $P( |\xi
- \eta| > x)$  for large thresholds $x$.
 {We gain more precision from HRV.} 

Define the cone 
$$\bO_1 = [0,\infty)^2 \setminus \bigl(\{\bzero\} \cup \text{DIAG}\bigr)
=\left\{ (u, v) \in
  \mathbb{D}: |u-v| > 0 \right\}.$$ 
The distribution of ${\bZ }$ has HRV on $\bO_1$ with scaling
function $b_1(t) = t^{1/2}$ and limit measure 
$$\nu_1(du,dv)= u^{-3}\epsilon_{((0,\infty)\times \{0\})} (du,
dv) 
+ u^{-3}\epsilon_{\bigl(\{0\}\times (0,\infty)\bigr) } (du, dv) ,\quad
(u,v)\in \bO_1,
$$ {the measure that concentrates mass on the axes and which is
restricted to $\bO_1$. This measure results from the second summand in
$\bZ$, namely $(1-B)(X_2,X_3)$; the first summand $B(X_1^2,X_1^2)$
contributes nothing to the limit due to the restiction to $\bO_1$.} So for instance, for $(u,v) \in \bO_1,\,x>0 ,$
$$\nu_1 \left( [{\bzero }, (u,v)]^c \cap \left \{ (u, v) \in
    \bO_1: |u - v| > x \right \} \right) = \frac1{2} \left(
  {\left(u \vee x\right)}^{-2} + {\left(v \vee x\right)}^{-2}
\right)$$ 
and  for some large $t > 0$, letting $u\downarrow 0$, $v\downarrow 0$,
we see
$$P ( |\xi- \eta| > x) \approx \frac1{t}{\left(\frac{x}{b_1(t)} \right)}^{-2}.$$
{Statistical estimates of the risk region probability
replace}  $b_1(t)$  by a statistic as in Section
\ref{sec:detection_estimation}.

The fact that $\nu_1(\cdot)$  concentrates on the
axes suggests seeking a further HRV property on a  cone smaller than $\bO_1$.
If one is needs risk probabilities of the form $P( \xi - \eta>
x, \eta > y)$ for large thresholds $x$ and $y$, we seek 
  HRV property, say, $\mathbb{O}^2 = \left\{ (u,v) \in \bO_1:
u,v> 0 \right\}$ or a sub-cone of $\mathbb{O}^2$. 

As an example of why  risk probabilities like $P(\xi -\eta > x)$
arise,  imagine investing in financial instruments $I_1$ and $I_2$
that  have risks $\xi$ and $\eta$ per unit of investment
 and suppose these risks have asymptotic full dependence.

For any $a_1,
a_2, c > 0$, asymptotic full dependence of $\xi$ and $\eta $ implies
$P(a_1 \xi + a_1 \eta > x)$ {should} be bigger
than $P(c(\xi - \eta) > x)$, provided  $x$ {is large}. So, if $l_1 >
l_2$,  it is less risky to invest in 
the financial instrument $I_1 - I_2$ rather than investing in both
$I_1$ and $I_2$. Obviously, investing in the financial instrument $I_1
- I_2$ requires us to measure risks associated with this portfolio,
which leads to the need  to evaluate $P(\xi - \eta > x)$ for large thresholds
$x$. 

In summary, 
this example shows how a more flexible definition of HRV possibly
allows computation of risk probabilities in the presence of asymptotic
full dependence.
}
\end{ex}

\begin{ex}\label{ex:infinite_sequence}
{\rm{ {\bf{Infinite sequence of cones:}} 
HRV was originally defined {for $d=2$} and for two cones \citep{resnick:2002a}
and then extended to a finite sequence of cones
\citep{mitra:resnick:2010}. In this paper our definition of HRV 
{allows} the possibility that we progressively find HRV  on an infinite sequence of cones.
{We exhibit an example for $d=2$ where this is indeed the case.}
An infinite sequence of cones may create problems for risk estimation
 which we discuss afterwards.

Suppose {$\{ X_i, i \geq 1\}$} are iid random variables
with common
Pareto($1$) distribution. Let {$\{ Y_1,Y_2 \}$} be iid with
common Pareto($2$) distribution  and suppose $\{ B_i, i \geq 1\}$ is an
infinite sequence of random variables with $P(B_i = 
1) = 1 - P(B_i = 0) = 2^{-i}$ and $\sum_{i = 1}^\infty B_i =
1$. {(For instance, let $T$ be the index of the first success in
  an iid sequence of Bernoulli trials and then set $B_i=1_{[T=i]}, i
  \geq 1.$)}
Assume that $\{ X_i, i\geq 1 \}$, $\{ Y_1,Y_2 \}$
and $\{ B_i : i\geq 1 \}$ are mutually independent. Define
the random vector ${\bZ }$ as 
$${\bZ } = (\xi, \eta) = B_1(Y_1, Y_2) + \sum_{i = 1}^\infty B_{i + 1}\left({(X_i)}^{\frac{1}{2-2^{-(i-1)}}},2^{i-1}{(X_i)}^{\frac{1}{2-2^{-(i-1)}}} \right). $$

So, ${\bZ }$ has regular variation on the cone $\bO=\bO_0 = \mathbb{D} =
[0, \infty)^d \backslash \{ {\bzero } \}$ with {index of regular
  variation $\alpha = 1$, scaling function $b(t)=t$,}
limit measure
$\nu(\cdot)$ {concentrating} on the diagonal $\text{DIAG}:=\{ (x,
x) : x \in [0, \infty)\}$. So, we remove the diagonal and {find} HRV on
$\mathbb{O}_{1} = \bD \backslash \text{DIAG}$ {with $b_1(t)=t^{2/3},\, \alpha_1=3/2$ and limit measure
$\nu_1$ concentrating on {the ray} $\{(x,2x): x\geq 0\}$.} 
 Progressively seeking HRV on 
successive cones, we find at the $(i 
+1)$-th step of our analysis, ${\bZ }$ has regular variation on the
cone  
\begin{equation}\label{eqn:oh}
\mathbb{O}_{i} = \bD \backslash \left[ \cup_{j=1}^i \{ (x, 2^{j-1}x) :
  x \in [0, \infty)\} \right]
\end{equation}
with the limit measure $\nu_i(\cdot)$ and the index of regular
variation $\alpha_i = 2 - 2^{-i}$. The limit measure $\nu_i (\cdot)$
concentrates on $\{(x, 2^ix) : x \in [0, \infty) \}$. 

{Selection of cones must be guided by the type of risk probability needed.
Consider trying to estimate} $P(\xi - \eta > x)$ for large thresholds
$x$ {using the cones $\bO_i, i\geq 0$ given in \eqref{eqn:oh}.}
At the $(i+1)$st stage{, {using cone $\bO_i$,}
the limit measure $\nu_i(\cdot)$ puts zero mass on the cone $\{ (u,
v): u > v \}$. So, even after a million HRV steps, 
we will estimate $P(\xi - \eta > x)$ for large thresholds
$x$ as zero, which is clearly wrong due to the definition
of ${\bZ }$ {since }
$$P[\xi-\eta>x]\geq P[B_1=1]P[Y_1-Y_2>x]>0.$$
An alternative procedure
seeks regular variation on
the cone $\{ (u, v): u > v \}$ and this {leads to somewhat}
more reasonable estimates of $P(\xi - \eta > x)$ for large thresholds
$x$ since in this case the regular variation with the Pareto(2)
variables is captured.

The moral of the story is that the choice of sequence of
cones when defining HRV should  be guided by the kind of risk
sets considered.
 For example, if we are only interested in
joint tail probabilities, a possible choice of sequence of cones is
$\mathbb{D} = \mathbb{D}^{(1)} \supset \mathbb{D}^{(2)} \supset \cdots
\supset \mathbb{D}^{(d)}$, where 
$$ \mathbb{D}^{(l)} = \{ \bx \in [0, \infty)^d : x^{(l)} > 0 \},
$$ 
and recall  $x^{(l)}$ is the $l$-th largest component of
  $\bx$.
This special case is discussed in \cite{mitra:resnick:2010}. 

}}}
\end{ex}

\begin{ex}\label{ex:CEV_id}
{\rm  {\bf The CEV model and mass on the lines through $\binfty$.} 
If we consider the CEV model on $\E_\sqcap$,
there can exist two
different limits in \eqref{eqn:CEV} under two different
normalizations. This problem disappears if we restrict convergence
to $\bD_\sqcap$.

Suppose $Y$ is Pareto(1)  and 
 $B$ is a  Bernoulli random variable with $P[B=1]=P[B=0]=1/2$. Define  
$$\bZ = (\xi,\eta) = B(Y,Y) + (1-B) (\sqrt{Y},Y).$$
Then the following convergences both hold in $\bE_{\sqcap}$: for $x\ge0,y>0$, 
\begin{align}
\nu_1([0,x]\times (y,\infty]):&=\lim\limits_{t \to \infty} tP\left[\left(\frac{\xi}{t},\frac{\eta}{t}\right) \in [0,x]\times (y,\infty]\right]\nonumber\\
                              &= \frac 12 \left(\frac1y -\frac 1x\right)_{+} +  \frac 1{2y}. \label{eq:ex:tt}\\
\nu_2([0,x]\times (y,\infty]):&=\lim\limits_{t \to \infty} tP\left[\left(\frac{\xi}{\sqrt{t}},\frac{\eta}{t}\right) \in [0,x]\times (y,\infty] \right]\nonumber\\
                              & = \frac 12 \left(\frac1y -\frac 1{x^2}\right)_{+} . \label{eq:ex:tthalf}
\end{align}
So $\bZ$ follows a CEV model on $\bE_{\sqcap}$ with two different scalings. Note that $\nu_1$ does not put any mass on lines through
$\binfty$, but  $\nu_2$ does: 
$ \nu_2(\{\infty\} \times(y,\infty]) = {1}/{2y}.$
If we restrict convergence to $\bD_\sqcap$, limits are unique
and $\bZ$ is regularly varying on $\bD_{\sqcap}$ with limit measure
$\nu_1$ given by \eqref{eq:ex:tt} restricted to $\bD_{\sqcap}$.

}
\end{ex}

\section{Estimating the spectral measure and its
  support}\label{sec:detection_estimation} 
We have defined regular variation on a big cone $\bO_0\subset\R^d$
along with hidden regular variation in a nested sequence of subcones
$\bO_0 \supset \bO_1 \supset \bO_2 \supset \ldots$. 
{We now propose strategies for deciding whether HRV is consistent
  with a given data set and, if so, how to} 
 estimate probabilities of sets pertaining to joint occurence of extreme or high values. 
{We proceed as follows:}
\begin{enumerate}
 \item Specify a fixed finite sequence of cones pertinent to
 the problem, and {seek}
HRV sequentially on these cones. We discuss this 
in Section \ref{subsec:spec_seq} which 
follows ideas proposed in \citep[Section 3]{mitra:resnick:2010}.
\item If the sequence of cones is not clear,
proceed by estimating the support of the limit measure at each
 step,
 removing it, {and seeking HRV on}
 the complement of the support.
Then the hidden limit measure is  estimated using semi-parametric techniques similar to \ref{subsec:spec_seq}. 
\end{enumerate}

\subsection{Specified sequence of cones}\label{subsec:spec_seq} 
Suppose $\bZ_1,\bZ_2,\ldots,\bZ_n$ are iid random vectors in $\bC \subset
[0,\infty)^d$ {whose common
distribution has a regularly varying tail on}   $\bO$ according to Definition \ref{def:closedconecomp}
 with normalizing function $b(\cdot)$
and limit measure $\nu(\cdot)$. Also assume that we have a 
specified sequence of cones $\bO=\bO_0 \supset \bO_1 \supset \bO_2
\supset \ldots$  where we seek regular variation. Such a sequence of cones is 
known and fixed.

We provide an estimate for the limit measure of regular
variation on $\bO$ and the same method can be applied to find limit
measures for hidden regular variation on the subcones.}

Now, according to Corollary \ref{cor:decomp2statements} regular
variation on $\bO$ as above is equivalent to assuming
{$P[d(\bZ,\bF)>x]$ is regularly varying at $\infty$} with some exponent $\alpha >0$ {and normalizing function $b(t)$ which
we take as $b(t):=F_d^{\leftarrow}(1-1/t)$ where $d(\bZ,\bF)$ has distribution function $F_d$} and 
\begin{equation}\label{eqn:whew}
 P\left[\frac{\bZ}{d(\bZ,\bF)}\in \;\cdot \; \Bigg| d(\bZ,\bF )> t
\right] \Rightarrow S_{\bO}(\cdot) \qquad (t\to \infty),
\end{equation}
in $P(\partial \aleph_{\bO})$, the class of all probability measures
on  $\partial \aleph_{\bO}= \{\bx\in \bO: d(\bx,\bF)=1\}$. Thus
 we estimate $\nu$ by estimating $\alpha$ and $S_{\bO}$
separately.
Considering $d_1^{\bF}:=d(\bZ_1,\bF), \ldots, d_n^{\bF}=d(\bZ_n,\bF)$
as iid samples from a regularly varying distribution on $(0,\infty)$,
the exponent $\alpha$ can be estimated 
using the Hill, Pickands or QQ estimators
\citep{resnickbook:2007}.  

{Here is an outline of how to obtain  an empirical estimator} of
$S_{\bO}$ following \citep{resnickbook:2007}. 
\begin{prop}\label{prop:EstSO}
Assume the common distribution of the iid random vectors $\bZ_1,
\dots, \bZ_n$ satisfies Definition \ref{def:closedconecomp} and \eqref{eqn:ccc}.
As $n \to \infty$, $k \to \infty$, $n/k \to \infty$, we have in  $P(\partial\aleph_{\bO})$,
\begin{align}\label{eqn:typeset}
S_n(\cdot):= \frac{\sum\limits_{i=1}^n
  \epsilon_{{\bigl(d_i^{\bF}}/b(n/k),{\bZ_i}/{d_i^{\bF}} \bigr)}
  ((1,\infty)\times\;\cdot\;)}{\sum\limits_{i=1}^n
  \epsilon_{{d_i^{\bF}}/b(n/k) }  (1,\infty)} \Rightarrow
S_{\bO}(\cdot). 
\end{align}
\end{prop}

\begin{proof}
Since $\{d_i^{\bF}, 1\le i \le n\}$ are iid 
regularly varying random variables from some distribution  $F_d$ on
$(0,\infty)$ with norming function $b(t)$ which can be chosen to be
$b(t)=F_d^{\leftarrow}(1-1/t)$, by Theorem
\ref{thm:relation_to_MRV}. Thus for $x>0$,
$\frac nk P[ d_i^{\bF} /b(n/k) >x]\to cx^{-\alpha},$
and from \citet[page 139]{resnickbook:2007}
 this is equivalent to 
$\frac 1k \sum_{i=1}^n \epsilon_{  d_i^{\bF}/b(n/k)        }
(1,\infty) \Rightarrow c,$
and to prove \eqref{eqn:typeset}, it suffices to show in  $M_+(\partial\aleph_{\bO})$,
\begin{equation}\label{eqn:sorta}
\frac 1k \sum_{i=1}^n \epsilon_{ \bZ_i/d_i^{\bF}       } (\cdot) 1_{[d_i^{\bF}/b(n/k)  >1]}
\Rightarrow 
S_{\bO}(\cdot). 
\end{equation}
The counting function on the left of \eqref{eqn:sorta} only counts
$\bZ_i/d_i^{\bF} $ such that $d_i^{\bF} /b(n/k)>1$. The distribution
of such random elements is 
$P[\bZ_i/d_i^{\bF}  \in \cdot | d_i^{\bF} /b(n/k)>1],$
\citep[page 212]{resnickbook:2008} and 
\eqref{eqn:whew} holds. Using \citet[Theorem 5.3ii, page
139]{resnickbook:2007} and the style of argument in \citet[page
213]{resnickbook:2008}, we get \eqref{eqn:sorta}.
\end{proof}

{\rm  The estimator $S_n$  of  $S_{\bO}$  in Proposition
  \ref{prop:EstSO}  relies on $b(t)$ which is typically unknown but
$b(n/k)$  can be estimated.
Order  $d_1^{\bF}, \ldots, d_n^{\bF}$  as  $d_{(1)}^{\bF}\ge  \ldots \ge d_{(n)}^{\bF}$
and  $d_{(k+1)}^{\bF}/b(n/k) \stackrel{P}{\to} 1$ so  $d_{(k+1)}^{\bF}$ is a consistent
estimator of $b(n/k)$ as $n\to \infty,k\to \infty, n/k\to
\infty$ \citep[page 81]{resnickbook:2007}. 
 Hence {we replace}
$b(n/k)$ by $d_{(k+1)}^{\bF}$ {and} propose the estimator $\hat{S}_n$ for $S_{\bO}$  in $S_n$  as follows:
\begin{align*}
 \hat{S}_n(\cdot):= \frac{\sum\limits_{i=1}^n
   \epsilon_{\{{d_i^{\bF}}/d_{(k+1)}^{\bF},{\bZ_i}/{d_i^{\bF}} \}}
   ((1,\infty)\times\;\cdot\;)}{\sum\limits_{i=1}^n
   \epsilon_{\{{d_i^{\bF}}/d_{(k+1)}^{\bF} \}}  (1,\infty)}
  = \frac 1k  \sum\limits_{i=1}^{{n}} {1_{[d_i^{\bF}/d_{(k+1)}^{\bF} >1]}}\epsilon_{{\bZ_i}/{d_i^{\bF}}}  (\cdot).
\end{align*}
}

\begin{prop}\label{prop:EstSOwithSnhat}
 As $n \to \infty, k \to , n/k \to \infty$, $\hat{S}_n \Rightarrow S_{\bO}$
 in $P(\partial\aleph_{\bO})$.
\end{prop}
\begin{proof}
The proof is a consequence of the continuous mapping theorem and
{Proposition \ref{prop:EstSO}. For  details  see, for instance, \citep{dehaan:resnick:1993}.}
\end{proof}
 
Thus  when $\bZ_1,\bZ_2,\ldots,\bZ_n$ are iid random vectors in
 $\bC \subset [0,\infty)^d$ which {have a  regularly varying}
 distribution on $\bO$, we can estimate both $\alpha$ and the spectral measure $S_{\bO}$. If
 we have a specified finite sequence of cones $\bO:=\bO_0 \supset
 \bO_1 \supset \bO_2 \supset \ldots \bO^m$, 
then we sequentially estimate  the limit measure by separately estimating the spectral measure and the index $\alpha_i$.

\subsection{Support estimation}\label{subsec:supp_est}
 Suppose $\bZ_1,\bZ_2,\ldots,\bZ_n$ are iid random vectors in $\bC
 \subset [0,\infty)^d$ {whose common distribution} is  regularly varying on $\bO$ according
 to Definition \ref{def:closedconecomp} 
with normalizing function $b(\cdot)$ and {limit} measure
$\nu(\cdot)$. 
{Without} a sequence of cones where hidden regular variation
can be sought, 
{the  task of exploring} for  appropriate cones
where HRV may exist {is challenging}. {One clear strategy is to identify the
  support of $\nu$, which we call $\text{supp}(\nu)$, and then seek HRV
  on the complement of the support.} Since
\begin{align}\label{eq:supp}
\supp (\nu) = \{t \cdot \supp (S_\bO ),\, t\geq 0\}.
\end{align}
it suffices to determine the support of $S_{\bO}$.

We
propose estimating  the support of the spectral measure $S_{\bO}$
with a point cloud; that is, a discrete random closed set.

\begin{prop}\label{prop:supportest} 
 Suppose $\bZ_1,\bZ_2,\ldots,\bZ_n$ are iid random vectors in $\bC
 \subset [0,\infty)^d$ whose common distribution is regularly varying on $\bO$ 
with normalizing function $b(\cdot)$ and limit measure $\nu(\cdot)$. 
As $n \to \infty, k \to \infty, n/k \to \infty$,
\begin{align}\label{convergence:supp}
 \supp_{k,n} & = \left\{\frac{\bZ_i}{d_i^{\bF}}: d_i^{\bF}>
   d_{(k+1)}^{\bF}, i =1,\ldots, n \right\} \Rightarrow \supp(S_{\bO}). 
\end{align}
\end{prop}

Convergence in \eqref{convergence:supp} occurs in the space of
closed sets under the Fell topology {or the space of compact sets in
the Hausdorff topology} \citep{molchanov:2005}.

\begin{proof} 
To show \eqref{convergence:supp}, it suffices from 
\cite[Proposition 6.10, page 87]{molchanov:2005} to show for any 
$h\in C_K^{+}(\partial\aleph_{\bO})$ that
\begin{equation}\label{eqn:diaper1}
E\Bigl(
 \sup_{i} \;\left\{ h\left(\frac{\bZ_i}{d_i^{\bF}}\right):
   d_i^{\bF}>d_{(k+1)}^{\bF}, 1\le i \le n\right\} \Bigr)
{\to}  \sup_{x} \;\{h(x): x \in \supp(S_{\bO})\}.
\end{equation}
From Proposition \ref{prop:EstSOwithSnhat},
$$
 \hat{S}_n(\cdot):= 
  \frac 1k  \sum\limits_{i=1}^{{n}} {1_{[d_i^{\bF}/d_{(k+1)}^{\bF}
      ]}}\epsilon_{{\bZ_i}/{d_i^{\bF}}} (\cdot).
\Rightarrow S_\bO,$$
and from the continuous mapping theorem, for any 
$h\in C_K^{+}(\partial\aleph_{\bO})$, we get in $P(\mathbb{R})$, the class of probability measures on $\R$,
\begin{equation}\label{eqn:diaper2}
 \hat{S}_n \circ h^{-1} =   \frac 1k  \sum\limits_{i=1}^{{n}} {1_{[d_i^{\bF}/d_{(k+1)}^{\bF}
      ]}}\epsilon_{{h(\bZ_i}/{d_i^{\bF}})} (\cdot)
\Rightarrow S_\bO \circ h^{-1}.
\end{equation}

If $F_n, n \ge 0$ are probability measures on
$\R$ with bounded support and $F_n \Rightarrow F_0$ then
\begin{align}\label{eq:convofrtendpt}
 x_{F_n}:=\sup_{x} \;\{x: F_n(x)<1\} & \to \sup_{x} \;\{x: F(x)<1\}=:x_F.
\end{align}
Applying this remark to \eqref{eqn:diaper2} and using the continuous
mapping theorem yields
 as $n\to \infty, k \to \infty, n/k \to \infty$, 
\begin{align}\label{eq:convofrtendpt_of_h}
 \sup_{i} \;\left\{ h\left(\frac{\bZ_i}{d_i^{\bF}}\right):
   d_i^{\bF}>d_{(k+1)}^{\bF}, 1\le i \le n\right\} & 
\Rightarrow  \sup_{x} \;\{h(x): x \in \supp(S_{\bO})\}.
\end{align}
Since $h\in C_K^{+}(\partial\aleph_{\bO})$ is always bounded above,
use dominated convergence applied to convergence in distribution to
get the desired \eqref{eqn:diaper1}.
\end{proof}

 Proposition \ref{prop:supportest} provides an estimate of
 $\supp(S_{\bO})$ and hence $\supp(\nu)$. In principle, we can remove
 the estimated support
 from $\bO$ and look for hidden regular variation in the complement.
How well this works in practice remains to be seen. For one thing, the
estimated support set of $S_\bO$ is always discrete meaning that the
estimated support of $\nu$ is a finite set of rays.  With a large data
set, we might be able to get  a fair idea about the support  
of the distribution and where to look for further hidden regular
variation. If there were reason to believe or hope that the support of
$S_\bO$ were convex, our estimation procedure could be modified by
taking the convex hull of the points in \eqref{convergence:supp}.

\section{Conclusion} \label{sec:conclusion}
Our treatment of regular variation on cones which is determined by the support
of the limit measures unifies under one theoretical  umbrella several
related concepts: asymptotic independence, asymptotic full dependence,
and the conditional extreme value model. Our approach highlights the
structural similarities of these concepts while making plain in what
ways the cases differ. Furthermore, the notion of $\bM^*$-convergence introduced
in Section \ref{sec:Mstar} provides a tool to deal with the generalized notion of regular
variation given here. Generalizing this notion 
of convergence and analyzing its properties admits potential for further research.

It is always an ambitious undertaking to
statistically identify lower 
order behavior and this project has not attempted data analysis or
{tested} the feasibility of the statistical methods
discussed in Section \ref{sec:detection_estimation}.
It is clear further work is required, particularly
for the case where the support of the limit measures must be
identified from data. One can imagine that for high dimensional data
whose dimension is of the order of hundreds,  sophistication is
required to pursue successive cones where regular variation exists.

\section{Appendix}\label{sec:appendix}
\subsection{Regular variation on $\E:=[0,\infty]^d\setminus
  \{\bzero\}$ vs $\bD=[0,\infty)^d
\setminus \{\bzero\}$}\label{subsec:equivRegVar}
We verify that the traditional notion of multivariate regular variation given in
  \eqref{eqn:existing_reg_var} on $\E$ is equivalent to Definition
  \ref{def:closedconecomp} if we choose $\bC = [0, \infty)^d$ and
  $\bF=\{\bzero\}$. 
This yields $\bO = \bD$.

\begin{thm}\label{thm:relation_to_MRV}
{{ Regular variation on $\mathbb{D}$ according to Definition \ref{def:closedconecomp} 
is equivalent to the {traditional} notion of multivariate regular variation
given in \eqref{eqn:existing_reg_var} and the limit measures
$\nu (\cdot)$ of \eqref{eqn:ccc} and $\mu(\cdot)$ of
\eqref{eqn:existing_reg_var} are equal on $\mathbb{D}$. {Moreover,
  $\mu(\cdot)$ puts zero measure on $\bE \setminus \bD$.} 
}}
\end{thm}

\begin{proof}
First we show that the standard notion of multivariate regular
variation {on $\E$} given in \eqref{eqn:existing_reg_var} implies 
\eqref{eqn:ccc} in $\bM^*(\bD)$.
 Let $\nu(\cdot)$ be a measure on $\mathbb{D}$ such that $\nu(\cdot) = \mu(\cdot)$. From \citet[page 176]{resnickbook:2007}, we get that $\mu(\E \backslash \mathbb{D}) = 0$. So, since $\mu(\cdot) \ne 0$ and non-degenerate, $\nu(\cdot) \ne 0$ and non-degenerate.

{For $B\subset \bD$,} note that $\partial B = \bar B
\backslash B^o$ is defined with respect to the relative topology on
$\mathbb{D}$ and hence $\partial B \subset \mathbb{D}$. Thus,
$\nu(\partial B) = 0$ implies $\mu( \partial B) = 0$.  Also, since
$[0, \infty]^d$ is a compact space, any set $B \subset \mathbb{D}$
bounded away from $\{ {\bzero } \}$ is a relatively compact set in $\E$
\cite[page 171, Proposition 6.1]{resnickbook:2007}. Therefore, by
definition, $\nu(\cdot) \in \bM^*([0,\infty)^d, \bD)$ and by
\eqref{eqn:existing_reg_var}, for any $B \subset \mathbb{D}$ bounded
away from $\{ {\bzero }\}$ and $\nu(\partial B) = 0$, 
$$tP\left[ \frac{{\bZ }}{b(t)} \in B \right] \rightarrow \mu(B) = \nu(B).$$
So, \eqref{eqn:ccc} holds with $\bC = [0, \infty)^d$, $\bO = \bD$, 
and {$b$ is the same as in \eqref{eqn:existing_reg_var} and $\nu $ is
the restriction of $\mu$ to $\bD$.}

{Conversely,} we show that 
Definition {\ref{def:closedconecomp} and \eqref{eqn:ccc}}
 {with $\bO = \bD$} implies the traditional  notion of multivariate
 regular variation on $\E$ in \eqref{eqn:existing_reg_var}. 
Define a measure {$\mu(\cdot)$ on $\E$ as $\mu(\cdot) = \nu(\cdot \cap
   \mathbb{D})$.}  A relatively compact set $B$ of
 $\E$ must be bounded away from $\{ {\bzero }\}$ \cite[page 171,
 Proposition 6.1]{resnickbook:2007}. So, from definition of
 $\mu(\cdot)$, it is Radon. Note that $\partial B = \bar B
 \backslash B^o$ is defined with respect to the topology on
 $\mathbb{E}$, but $\partial (B \cap \mathbb{D})$ is defined with
 respect to the relative topology on $\mathbb{D}$. Also from the
 definition of $\mu(\cdot)$, $\mu(\partial B) = 0$ implies
 $\nu(\partial (B \cap \mathbb{D})) = \nu( \partial B \cap \mathbb{D})
 = \mu(\partial B) = 0$. Therefore, from \eqref{eqn:ccc},
 for any relatively compact set $B$ of $\E$
 such that $\mu(\partial B) = 0$, as $t \rightarrow \infty$, 
$$tP\left[ \frac{{\bZ }}{b(t)} \in B \right] = tP\left[
  \frac{{\bZ }}{b(t)} \in B \cap \mathbb{D} \right] \rightarrow
\nu(B \cap \mathbb{D}) = \mu(B).$$ 
The first equality above holds since ${\bZ } \in [0,
\infty)^d$. Hence, vague convergence in \eqref{eqn:existing_reg_var}
holds {with 
the same $b$ as in \eqref{eqn:ccc} and with $\mu$ as the extension of
$\nu$ from $\bD$ to $\E$.}
\end{proof}

Regular variation on $\bD$ can also be expressed in terms of the polar
coordinate transformation. As at the begining of Section
\ref{subsec:equivRegVar}, set 
 $\bC = [0, \infty)^d$,
  $\bF=\{\bzero\}$ and  $\bO = \bD$.

\begin{prop}\label{prop:bigconedecomp} 
Regular variation on $\mathbb{O}$ as given in Definition \ref{def:closedconecomp} is equivalent to the following condition:
\begin{equation}\label{eqn:bigcone_decomp}
tP\left[ \left( \frac{||{\bZ }||}{b(t)},
    \frac{{\bZ }}{||{\bZ }||} \right) \in A \right] \rightarrow
\nu_{(\alpha )} \times S_{\mathbb{O}}(A) 
\end{equation}
for all measurable $A \subset (0, \infty) \times \partial
\aleph_{\mathbb{O}}$ such that $\overline{A^1} \cap \{ 0 \} =
\emptyset$ and $\nu \circ h^{-1} (\partial A) = 0$, where $A^1$ is the
projection of $A$ on its first coordinate, $h(\cdot)$ is a function
defined by $h : \bx \mapsto \left(||\bx||,
  \frac{\bx}{||\bx||} \right)$, $\partial
\aleph_{\mathbb{O}} = \{\bx \in \mathbb{O}: ||\bx|| = 1\}$,
$S_{\mathbb{O}}(\cdot)$ is a probability measure on $\partial
\aleph_{\mathbb{O}}$ and $\nu_{(\alpha )}(\cdot)$ is a Pareto measure
given by $\nu_{(\alpha )}((x, \infty) ) = x^{-\alpha}$ for $x > 0$. The probability measure $S_{\mathbb{O}}(\cdot)$ is called the spectral measure and is related to $\nu(\cdot)$ by the relation
\begin{equation}\label{eqn:sandnu}
S_{\mathbb{O}}(\Lambda) = \nu \left( \left\{ \bx \in \mathbb{O}: ||\bx|| \ge 1, \frac{\bx}{|| \bx ||} \in \Lambda \right\} \right).
\end{equation}
\end{prop}

\begin{proof} {This is a special case of}    Proposition \ref{prop:smallconedecomp}. \end{proof}

\subsection{Regular variation on $\E^{(l)}=[0,\infty]^d\setminus \{\bx
  \in \E: x^{(l)}>0\}$ vs $\bD^{(l)}=[0,\infty)^d\setminus \{\bx
  \in [0,\infty)^d:x^{(l)}>0\}$} \label{rem:relation_with_existing}
Recall $x^{(l)}$ is the $l$-th largest component of $\bx$,
 $l = 1, 2, \cdots, d$. Hidden regular variation using $
\E^{(l)}$ is considered in   \citet{mitra:resnick:2010}. Unlike the
situation in subsection \ref{subsec:equivRegVar}, here limit measures
can put mass on $\E^{(l)} \setminus \bD^{(l)} $ as found in 
 \citet{mitra:resnick:2010}. We compare regular variation in
 $\E^{(l)}$ using the traditional vague convergence definition in
 which the vague convergence in \eqref{eqn:existing_reg_var} is
 assumed to hold in $\bM_+(\E^{(l)})$ with regular variation given in
 \eqref{eqn:ccc} in $\bM^*(\bC,\bO)$ where
$\bC=[0,\infty)^d, $ $\bF=\{\bx \in [0,\infty)^d: x^{(l)}=0\} $ and 
$\bO=\bD^{(l)}=\bC\setminus \bF$.

\begin{thm}\label{thm:relation_to_RV_smallcone1}
Regular variation on $\bM^*(\mathbb{D}^{(l)})$ is equivalent to
the traditional vague convergence notion of regular variation in
$\bM_+(\E^{(l)})$
 if the limit measure
$\mu(\cdot)$ given in the $\bM_+(\E^{(l)})$ analogue of
 \eqref{eqn:existing_reg_var} does not
give any mass to the set $\E^{(l)}\setminus \bD^{(l)}$. 
In this case,
the limit measures $\nu(\cdot)$ of \eqref{eqn:ccc} and
$\mu(\cdot)$ of 
\eqref{eqn:existing_reg_var} 
are equal on
$\mathbb{D}^{(l)}$. 
\end{thm}

\begin{proof}
Suppose, for a random vector ${\bZ }$, there exist
 a function $b^{(l)}(t) \uparrow \infty$  and 
 a non-negative non-degenerate Radon measure $\mu^{(l)}(\cdot) \ne 0$ on $\E^{(l)}$,
 such that in $\bM_+(\E^{(l)})$
\begin{equation}\label{eqn:existing_hidden_defn}
tP\left[ \frac{{\bZ }}{b^{(l)}(t)} \in \cdot \right] \stackrel{v}{\rightarrow} \mu^{(l)}(\cdot).
\end{equation}
and the limit measure
$\mu^{(l)}(\cdot)$ does not give any mass to $\E^{(l)} \setminus
\bD^{(l)}$. Define a measure $\chi(\cdot)$ on $\mathbb{D}^{(l)}$ as
$\chi(\cdot) = \mu^{(l)}(\cdot)$. Since $\mu^{(l)}(\cdot) \ne 0$ is
non-negative, non-degenerate and $\mu^{(l)}(\E^{(l)} \setminus
\bD^{(l)}) = 0$, the measure $\chi(\cdot) \ne 0$ is non-negative and
non-degenerate.  The subsets of $\mathbb{D}^{(l)}$ bounded
away from ${ \left( \mathbb{D}^{(l)} \right)}^c = [0, \infty)^d
\setminus \mathbb{D}^{(l)}$ are relatively compact in
$\E^{(l)}$. Therefore, using the fact that $\mu^{(l)}(\cdot)$ is Radon
and the definition of $\chi(\cdot)$, it follows that $\chi(\cdot)$
gives finite measure to sets bounded away from ${ \left(
    \mathbb{D}^{(l)} \right)}^c$. 
From the definition of
$\bM^*$-convergence, it follows that ${\bZ }$ satisfies
\eqref{eqn:ccc}
 with the scaling function $b(\cdot) =
b^{(l)}(\cdot)$ and the limit measure $\nu (\cdot) = \chi(\cdot)$. 

Conversely, suppose a random vector ${\bZ }$ satisfies
\eqref{eqn:ccc}  in $\bM^*(\bC,\bO)$
with $\bC = [0, \infty)^d$ and $\bO =
\bD^{(l)}$. Define a measure $\mu (\cdot)$ on $\E^{(l)}$ as
$\mu(\cdot) = \nu (\cdot \, \cap \, \mathbb{D}^{(l)})$. Since
$\nu(\cdot) \ne 0$ and is non-negative and non-degenerate, so is
$\mu(\cdot)$. A subset of $\E^{(l)}$ is relatively compact
in $\E^{(l)}$ iff it is bounded away from $\{ \bx \in
[0, \infty]^d: x^{(l)} = 0 \}$. Since $\nu(\cdot)$ gives
finite mass to sets bounded away from ${\left( \mathbb{D}^1
  \right)}^c$, from the definition of $\mu(\cdot)$, it follows that
$\mu(\cdot)$ is a Radon measure. From the description of the
compact sets in $\E^{(l)}$, it follows that ${\bZ }$ also satisfies
\eqref{eqn:existing_hidden_defn} with the scaling function
$b^{(l)}(\cdot) = b(\cdot)$ and the limit measure $\mu^{(l)}(\cdot)
= \mu(\cdot)$ \citep[page 52, Theorem 3.2]{resnickbook:2007}. 
\end{proof}

The set $\E^{(l)} \setminus \bD^{(l)} =
  \{ \bx \in \E^{(l)}: ||\bx||= \infty \}$ is the union of the lines through
  ${\boldsymbol{\infty}}$.  We emphasize that there 
exist examples of random vectors ${\bZ }$ which satisfy
\eqref{eqn:existing_hidden_defn} and the limit measure
{$\mu^{(l)}(\cdot)$ gives positive measure on the set $\E^{(l)}
  \setminus \bD^{(l)}$ \citep{mitra:resnick:2010}. 

\subsection{Regular variation on $\E_\sqcap=[0,\infty]\times
  (0,\infty)$ vs $\bD_\sqcap =[0,\infty)\times
  (0,\infty)$} \hfill

Recall CEV model from Section \ref{subsec:cev}.
\begin{prop}\label{prop:CEVHRV}
The following are equivalent:
\begin{enumerate}
 \item[(i)] $\bZ\in CEV(b_1,b_2)$ with limit measure $\mu(\cdot)$ and $b_1\sim b_2$ with 
\begin{align} 
& \mu\left([0,\infty]\times\{\binfty\} \cup \{\binfty\}\times(0,\infty]\right)=0.  \label{eqn:muinfis0}
\end{align}
\item[(ii)] $\bZ$ is regularly varying on $\bD_{\sqcap}$ according to
  \eqref{eqn:ccc}
 with  normalizing function $b_1$
and limit measure $\nu$ which does not concentrate on $\{0\}\times(0,\infty)$. 
\end{enumerate}
Also, if either of (i) or (ii) holds then $\mu(\cdot)=\nu(\cdot)$ on $\bD_{\sqcap}$. 
\end{prop}

\begin{proof}
 (i) implies (ii): Since $b_1\sim b_2$, \eqref{eqn:CEV} implies $\bZ
 \in CEV(b_1,b_1)$.
 Now, \eqref{eqn:CEV} implies
that for all relatively compact Borel sets $B$ in $\bD_\sqcap \subset \bE_{\sqcap}$ with $\mu(\partial B)=0$, 
\[tP\left[\frac{\bZ}{b_1(t)} \in B \right] \to \mu(B)\] 
as $t\to \infty$. Clearly $B$ is bounded away from $\bF$. Also $\mu$ is non-null and satisfies \eqref{eqn:muinfis0}. Thus $\nu(\cdot)=\mu(\cdot)|_{\bD_{\sqcap}}$ is
non-negative and non-degenerate on $\bD_{\sqcap}$. Hence $\bZ$ is regularly varying on $\bD_{\sqcap}$ with limit measure $\nu$. 
The non-degeneracy condition (a) for the CEV model implies that $\mu$
cannot concentrate on  $\{0\}\times(0,\infty)$. 
Conversely, if  (ii) implies (i) extend $\nu$ to a measure $\mu$ on $\bE_{\sqcap}$
which satisfies \eqref{eqn:muinfis0}.\end{proof}

\begin{rem}
{\rm We can drop the condition that $\mu $ does not concentrate on $\{0\}\times(0,\infty)$ in statement (ii) of Proposition \ref{prop:CEVHRV},
if we drop condition (a) from Definition \ref{CEV_defn} of the CEV model. }
\end{rem}

\bibliographystyle{plainnat}
\bibliography{/hg/u/bikram/Documents/Texfiles/bibfilenew}

 \end{document}